\documentclass[11pt]{report}
\usepackage{amsmath}
\usepackage{amssymb}
\usepackage{amsthm}
\usepackage{latexsym}
\usepackage{graphicx}
\usepackage{hyperref}

\newtheorem{thm}{Theorem}[chapter]
\newtheorem{cor}[thm]{Corollary}
\newtheorem{lem}[thm]{Lemma}
\newtheorem{prop}[thm]{Proposition}

\providecommand{\norm}[1]{\left\| #1 \right\|}

\newcommand{\mb}{\mathbf}
\newcommand{\mh}{\mathbb}
\newcommand{\mr}{\mathrm}
\newcommand{\mc}{\mathcal}
\newcommand{\mf}{\mathfrak}

\newcommand{\scs}{\scriptstyle}
\newcommand{\ts}{\textstyle}

\newcommand{\isom}{\xrightarrow{\;\sim\;}}

\newcommand{\cc}{\mathrm{cc}}
\newcommand{\es}{\emptyset}

\newcommand{\al}{_\mathrm{alg}}
\newcommand{\sm}{_\mathrm{smooth}}
\newcommand{\pch}{periodic cyclic homology }
\newcommand{\hot}{\widehat{\otimes}}

\newcommand{\xnr}{X_{\mathrm{nr}}}
\newcommand{\xunr}{X_{\mathrm{unr}}}
\newcommand{\inp}[2]{\langle #1 \,,\, #2 \rangle}
\newcommand{\prefix}[3]{\vphantom{#3}#1#2#3}
\newcommand{\matje}[4]{\left(\begin{smallmatrix} #1 & #2 \\ 
#3 & #4 \end{smallmatrix}\right)}
\newcommand{\hexagon}[6]{\begin{array}{ccccc}
  #1 & \to & #2 & \to & #3 \\
  \uparrow & & & & \downarrow \\
  #4 & \leftarrow & #5 & \leftarrow & #6
  \end{array}}

\begin{document}

\ \vspace{1cm}
\begin{center}
\huge \textbf{Periodic cyclic homology\\ of reductive $p$-adic groups}\\[1cm]
\Large Maarten Solleveld\\[3mm]
\normalsize Mathematisches Institut,
Georg-August-Universit\"at G\"ottingen\\
Bunsenstra\ss e 3-5, 37073 G\"ottingen, Germany\\
email: Maarten.Solleveld@mathematik.uni-goettingen.de \\[5mm]
August 2008 \\[3cm]
\end{center}

\begin{minipage}{13cm}
\textbf{Abstract.}\\
Let $G$ be a reductive $p$-adic group, $\mc H (G)$ its Hecke algebra
and $\mc S (G)$ its Schwartz algebra. We will show that these algebras
have the same periodic cyclic homology. This provides an
alternative proof of the Baum--Connes conjecture for $G$, modulo torsion.

As preparation for our main theorem we prove two results that have
independent interest. Firstly, a general comparison theorem for the 
periodic cyclic homology of finite type algebras and
certain Fr\'echet completions thereof. Secondly, a refined form of the 
Langlands classification for $G$, which clarifies the
relation between the smooth spectrum and the tempered spectrum.
\\[1cm]
\textbf{Acknowledgements.}\\
This paper is partly based on the author's PhD-thesis, which was written
at the Universiteit van Amsterdam under the supervision of Eric Opdam. 
The author is grateful for the support and advice that professor Opdam
has given him during his PhD-research.
He would also like to thank Ralf Meyer, Christian Voigt and the referee
for their comments, which lead to substantial clarifications of some proofs.
\end{minipage}

\tableofcontents

\chapter*{Introduction}
\addcontentsline{toc}{chapter}{Introduction}

In this paper we compare different homological invariants of algebras 
associated to reductive $p$-adic groups. Group algebras, or more precisely 
convolution algebras of functions on groups, have always been important 
objects of study in noncommutative geometry. Generally speaking the idea
(or hope) is that the interaction between representation theory, harmonic
analysis, operator algebras and geometry leads to results that can not (yet)
be proven inside only one of these areas. 

By definition a group algebra encodes information about a group, so its 
homological invariants should reflect properties of the group. 
Therefore, whenever one considers two convolution algebras 
associated to the same group, their invariants should be closely related. 
Yet in practice this has to be taken with quite a few grains of salt.
For example the \pch of $\mh C [\mh Z \rtimes C_2 ]$ is isomorphic
to the De Rham-cohomology (with complex coefficients) of the disjoint union of 
$\mh C^\times / (z \sim z^{-1})$ and a point. On the other hand the \pch of
the group-$C^*$-algebra $C^* (\mh Z \rtimes C_2 )$ does not give any new
information: it is the algebra itself in even degrees and it vanishes in odd 
degrees. So finding a meaningful invariant of the group is a matter of both
choosing the right group algebra and the right functor.

For Fr\'echet algebras topological $K$-theory is a good choice, since it is a 
very stable functor. It has the excision property and is invariant under homotopy 
equivalences and under passing to holomorphically closed dense subalgebras. 
Comparing with the above example, the $K$-theory of $C^* (\mh Z \rtimes C_2)$
is again isomorphic to the cohomology of a manifold. But the manifold
has been adjusted to its compact form 
\[
S^1 / (z \sim z^{-1}) \cup \mr{point} \quad \cong \quad [-1,1] \cup \mr{point}
\]
and we must take its singular cohomology with integral coefficients. We remark 
that subalgebras consisting of all functions on $\mh Z \rtimes C_2$ with rapid 
(resp. subexponential) decay have the same $K$-theory.

Nevertheless it can be hard to compute a $K$-group of a lesser-known algebra.
Indeed in the classical picture of $K_0$ one has to find all homotopy classes of 
projectors, a task for which no general procedure exists.
\\[1mm]

Of course there is a wider choice of interesting functors. Arguably the most subtle
one is Hochschild homology ($HH_*$), the oldest homology theory for algebras.
Depending on the circumstances it can be regarded as group cohomology,
(noncommutative) differential forms or as a torsion functor. Moreover Hochschild
homology can be computed with the very explicit bar complex. On the other hand
$HH_*$ does neither have the excision property, nor is it homotopy invariant. 

We mainly discuss \pch $(HP_* )$ in this paper. Although it carries less information 
than Hochschild homology, it is much more stable. The relation between $HH_*$ 
and $HP_*$ is analogous to that between differential forms and De Rham 
cohomology, as the Hochschild--Kostant--Rosenberg theorem makes explicit in the 
case of smooth commutative algebras. It is known that \pch has the excision 
property and is invariant under Morita equivalences, diffeotopy equivalences and
nilpotent extensions.
Together with the link to Hochschild homology these make $HP_*$ computable
in many cases. This functor works especially well on the category of finite type
algebras \cite{KNS}, that is, algebras that are finitely generated modules over the
coordinate ring of some complex affine variety. In this category an important 
principle holds for periodic cyclic homology, namely that it depends only on 
the primitive ideal spectrum of the algebra in question.

A similar principle fails miserably for topological algebras, even for commutative 
ones. For example let $M$ be a compact smooth manifold. Then 
$HP_* (C^\infty (M))$ is the De Rham cohomology of $M$, while $HP_* (C(M))$
just returns the $C^*$-algebra $C(M)$. The underlying reason is that $HP_*$
does not only see the (irreducible) modules of an algebra, it also takes the 
derived category into account. In geometric terms this means that $HP_* (A)$
does not only depend on the primitive ideal spectrum of $A$ as a topological 
space, but also on the structure of the "infinitesimal neighborhoods" of points
in this space. These infinitesimal neighborhoods are automatically right for 
finite type algebras, because they can be derived from the underlying affine variety.
But the spectrum of $C(M)$ does not admit infinitesimal neighborhoods. Indeed,
these have to be related to the powers of a maximal ideal $I$, but they collapse
because $\overline{I^n} = I$ for all $n \in \mh N$.

We remark that this problem can partially be overcome with a clever variation
on $HP_*$, local cyclic homology \cite{Mey3}. This functor gives nice results for
$C^*$-algebras because it is stable under isoradial homomorphisms of complete 
bornological algebras. On the other hand this theory does require an array of new 
techniques.

We will add a new move under which \pch is invariant. Let $\Gamma$ be a finite 
group acting (by $\alpha$) on a nonsingular 
complex affine variety $X$, and suppose that we have a cocycle 
$u: \Gamma \to GL_N (\mc O (X))$. Then $\alpha$ and $u$ combine to 
an action of $\Gamma$ on $M_N (\mc O (X))$:
\begin{equation}
\gamma \cdot f = u_\gamma f^{\alpha (\gamma)} u^{-1}_\gamma \,.
\end{equation}
The algebra of $\Gamma$-invariants $M_N (\mc O (X) )^\Gamma$ has a natural 
Fr\'echet completion, namely $M_N (C^\infty (X) )^\Gamma$. We will show in
Chapter 1 that the inclusion map induces an isomorphism
\begin{equation}\label{eq:0.1}
HP_* \big( M_N (\mc O (X) )^\Gamma \big) \to 
HP_* \big( M_N (C^\infty (X) )^\Gamma \big) \,.
\end{equation}
The proof is based on abelian filtrations of both algebras, that is, on sequences of ideals
such that the successive quotients are Morita equivalent to commutative algebras.
In terms of primitive ideal spectra this means that we have stratifications of finite length
such that all the strata are Hausdorff spaces.
\\[1mm]

Let us discuss these general issues in connection with reductive $p$-adic groups. 
We use this term as an abbreviation of ``the $\mh F$-rational points of a 
connected reductive algebraic group, where $\mh F$ is a non-Archimedean local 
field". Such groups are
important in number theory, especially in relation with the Langlands program.
There are many open problems for reductive $p$-adic groups, for example there is
no definite classification of irreducible smooth representations. There are two general 
strategies to divide the classification problem into pieces, thereby reducing it to 
either supercuspidal or square-integrable representations. 

For the first we start with a 
supercuspidal representation of a Levi-component of a parabolic subgroup of 
our given group $G$. Then we apply parabolic induction to obtain a (not necessarily 
irreducible) smooth $G$-representation. The collection of representations obtained
in this way contains every irreducible object at least once.

The second method involves the Langlands classification, which reduces the problem
to the classification of irreducible tempered $G$-representations. These can be found
as in the first method, replacing supercuspidal by square-integrable representations. 
This kind of induction was studied in \cite{ScZi}. The procedure yields a collection of 
(possibly decomposable) tempered $G$-representations,
in which every irreducible tempered representation appears at least once. 

Our efforts in Chapter 2 result in a refinement of the Langlands classification. To every 
irreducible smooth $G$-representation we associate a quadruple $(P,A,\omega, \chi )$
consisting of a parabolic pair $(P,A)$, a square-integrable representation $\omega$
of the Levi component $Z_G (A)$ and an unramified character $\chi$ of $Z_G (A)$.
Moreover we prove that this quadruple is unique up to $G$-conjugacy. This result is 
useful for comparing the smooth spectrum of $G$ with its
tempered spectrum, and for constructing stratifications of these spectra.
\\[1mm]

Let us consider three convolution algebras associated to a reductive $p$-adic group $G$.
Firstly the reduced $C^*$-algebra $C_r^* (G)$, secondly the Hecke algebra $\mc H (G)$
and thirdly Harish-Chandra's Schwartz algebra $\mc S (G)$. For each of these algebras 
we will study the most appropriate homology theory. For the reduced $C^*$-algebra this 
is topological $K$-theory, and for the Hecke algebra we take periodic cyclic homology. 
For the Schwartz algebra the choice is more difficult. Since it is not a Fr\'echet algebra 
the usual versions of $K$-theory are not even defined for $\mc S (G)$. It is not difficult 
to give an ad-hoc definition, and the natural ways to do so quickly lead to
$K_* (\mc S (G)) \cong K_* (C_r^* (G))$. Nevertheless, we would also like to compute
the \pch of $\mc S (G)$. It is definitely not a good idea to do this with respect to the 
algebraic tensor product, because that would ignore the topology on $\mc S (G)$. As
explained in \cite{Mey}, $\mc S (G)$ is best regarded as a bornological algebra, and
therefore we will study its \pch with respect to the completed bornological tensor product
$\hot_{\mh C}$. 

That this is the right choice is vindicated by two comparison theorems. 
On the one hand the author already proved in \cite{Sol1} that the Chern character 
for $\mc S (G)$ induces an isomorphism
\begin{equation}\label{eq:0.2}
ch \otimes \mr{id} : K_* (C_r ^* (G)) \otimes_{\mh Z} \mh C \to HP_* (\mc S (G) ,\hot_{\mh C}) \,.
\end{equation}
On the other hand we will show in Section \ref{sec:3.1} that the inclusion of $\mc H (G)$
in $\mc S (G)$ induces an isomorphism
\begin{equation}\label{eq:0.3}
HP_* (\mc H (G)) \to HP_* (\mc S (G) ,\hot_{\mh C}) \,.
\end{equation}
Of course both comparison theorems can be decomposed as direct sums over the
Bernstein components of $G$. The proof of \eqref{eq:0.3} is an extension of the ideas 
leading to \eqref{eq:0.1} and is related to the following quote \cite[p. 3]{SSZ}:\\
\emph{``The remarkable picture which emerges is that Bernstein's decomposition of
$\mc M (G)$ into its connected components refines into a stratification of $\mc G (G)$ 
where the strata, at least up to nilpotent elements, are module categories over 
commutative rings. We strongly believe that such a picture holds true for any group $G$."}\\
If this is indeed the case then our methods can be applied to many other groups.
\\[1mm]

The most important application of \eqref{eq:0.2} and \eqref{eq:0.3} lies in their relation
with yet other invariants of $G$. Let $\beta G$ be the affine Bruhat--Tits building of $G$.
The classical paper \cite{BCH} introduced among others the equivariant $K$-homology
$K_*^G (\beta G)$ and the cosheaf homology $CH_*^G (\beta G)$. Let us recall the
known relations between these invariants. The Baum--Connes conjecture for $G$,
proven by Lafforgue \cite{Laf}, asserts that the assembly map
\begin{equation}\label{eq:0.4}
\mu : K_*^G (\beta G) \to K_* (C_r^* (G))
\end{equation}
is an isomorphism. Voigt \cite{Voi4} constructed a Chern character 
\begin{equation}\label{eq:0.5}
ch : K_*^G (\beta G) \to CH_*^G (\beta G)
\end{equation}
which becomes an isomorphism after tensoring the left hand side with $\mh C$.
Furthermore it is already known from \cite{HiNi} that $CH_*^G (\beta G)$ is 
isomorphic to $HP_* (\mc H (G))$. Altogether we get a diagram
\begin{equation}\label{eq:0.6}
\begin{array}{ccc}
K_*^G (\beta G) \otimes_{\mh Z} \mh C & \cong & 
K_* (C_r^* (G)) \otimes_{\mh Z} \mh C \\
\cong & & \cong \\
CH_*^G (\beta G) & \cong & HP_* (\mc H (G))
\end{array}
\end{equation}
whose existence was already conjectured in \cite{BHP2}. We will prove in Section
\ref{sec:3.3} that it commutes. The four isomorphisms 
all have mutually independent proofs, so any three of them can be used to proof
the fourth. None of the proofs is easy, but it seems to the author that 
\eqref{eq:0.4} is the most difficult one. Therefore
it is not unreasonable to say that this diagram provides an alternative way to prove
the Baum--Connes conjecture for reductive $p$-adic groups, modulo torsion.

Returning to our initial broad point of view, we conclude that we used 
representation theory and harmonic analysis to prove results in noncommutative
geometry. It is outlined in \cite{BHP2} how cosheaf homology could be used to
prove representation theoretic results. The author hopes that the present paper
might contribute to the understanding of the issues raised in \cite{BHP2}.

\chapter{Comparison theorems for \pch}

\section{Finite type algebras}
\label{sec:1.1}

We will compare the \pch of certain finite type algebras and completions thereof. 
The motivating example of the result we aim at is as follows.

Let $X$ be a nonsingular complex affine variety. We consider the
algebras $\mc O (X)$ of regular (polynomial) functions and
$C^\infty (X)$ of complex valued smooth functions on $X$. By
default, if we talk about continuous or differentiable functions on
$X$ or about the cohomology of $X$, we always do this with respect
to the analytic topology on $X$, obtained from embedding $X$ in a
complex affine space.

The Hochschild--Kostant--Rosenberg--Connes theorem tells us what the
\pch of these algebras looks like:
\begin{equation}\label{eq:1.1}
\begin{array}{lll}
HP_n (\mc O (X)) & \cong &
\bigoplus_{m \in \mh Z} H_{DR}^{n+2m} (X ; \mh C) \,, \\
HP_n (C^\infty (X)) & \cong &
\bigoplus_{m \in \mh Z} H_{DR}^{n+2m} (X ; \mh C) \,.
\end{array}
\end{equation}
In the first line we must take the De Rham cohomology of $X$ as an
algebraic variety. However, according to a result of Grothendieck
and Deligne this is naturally isomorphic to the De Rham cohomology
of $X$ as a smooth manifold. Hence the inclusion $\mc O (X) \to
C^\infty (X)$ induces an isomorphism
\begin{equation}
HP_* (\mc O (X)) \to HP_* (C^\infty (X)) \,.
\end{equation}
Now let us discuss this in greater generality, allowing
noncommutative algebras. We denote the primitive ideal spectrum of
any algebra $A$ by Prim$(A)$ and we endow it with the Jacobson
topology, which is the natural noncommutative generalization of the
Zariski topology. An algebra homomorphism $\phi : A \to B$ is called
spectrum preserving if it induces a bijection on primitive ideal spaces,
in the following sense. For every $J \in \mr{Prim}(B)$ there is a 
unique $I \in \mr{Prim}(A)$ such that $\phi^{-1}(J) \subset I$,
and the map Prim$(B) \to \mr{Prim}(A) : J \mapsto I$ is bijective.

Since we do not want to get too far away from commutative
algebras, we will work with finite type algebras, see
\cite{KNS,BaNi}. Let $\mb k$ be the ring of regular functions on
some complex affine variety. A finite type $\mb k$-algebra is a
$\mb k$-algebra that is finitely generated as a $\mb k$-module.
The \pch of a finite type algebra always has finite dimension,
essentially because this is case for $\mc O (X)$ \cite[Theorem 1]{KNS}. 
Moreover it depends only on the primitive ideal spectrum
of the algebra, in the following sense:

\begin{thm}\label{thm:1.1} \textup{\cite[Theorem 8]{BaNi}}\\
A spectrum preserving morphism of finite type $\mb k$-algebras
induces an isomorphism on \pch \!.
\end{thm}

Morally speaking $HP_* (A)$ should correspond to the ``cohomology'' of Prim$ (A)$.
However, this is only a nonseparated scheme, so classical cohomology theories will 
not do. Yet this can be made precise with sheaf cohomology \cite[Section 2.2]{Sol3}.

It is not unreasonable to expect that there is always some
Fr\'echet completion $A\sm$ of $A = A\al$ such that the inclusion 
$A\al \to A\sm$ induces an isomorphism
\begin{equation}\label{eq:1.4}
HP_* (A\al ) \to HP_* (A\sm ) \,.
\end{equation}
A good candidate appears to be
\begin{equation}
A\sm = A\al \otimes_{\mc O (X)} C^\infty (X)
\end{equation}
if the center of $A\al$ is $\mc O (X)$. However I believe that it
would be rather cumbersome to determine precisely under which conditions this
works out. Moreover I do not know whether the resulting smooth
algebras are interesting in this generality. Therefore we restrict
our attention to algebras of a specific (but still rather general) form,
which we will now describe.

Let $\Gamma$ be a finite group acting (by $\alpha$) on the nonsingular
complex affine variety $X$. Take $N \in \mh N$ and consider the
algebra of matrix-valued regular functions on $X$:
\begin{equation}
\mc O (X ; M_N (\mh C) ) := M_N (\mc O (X)) =
\mc O (X) \otimes M_N (\mh C ) \,.
\end{equation}
Suppose that we have elements $u_\gamma \in GL_N (\mc O (X))$ such that
\begin{equation}\label{eq:1.2}
(\gamma \cdot f) (x) = u_\gamma (x) f (\alpha_\gamma^{-1} x)
u_\gamma^{-1} (x)
\end{equation}
defines a group action of $\Gamma$ on $M_N (\mc O (X))$, by
algebra homomorphisms. We do not require that $\gamma \mapsto
u_\gamma$ is a group homomorphism. Nevertheless the above does imply 
that there exists a 2-cocycle $\lambda : \Gamma \times \Gamma \to \mc O
(X)^\times$ such that
\[
u_\gamma (u_{\gamma'} \circ \alpha_\gamma^{-1}) =
\lambda (\gamma , \gamma') u_{\gamma \gamma'} \,.
\]
In particular, for every $x \in X$ we get a projective $\Gamma_x$-representation
\begin{equation}\label{eq:1.18}
(\pi_x ,\mh C^N ) \quad \mr{with} \quad \pi_x (\gamma ) = u_\gamma (x) \,.
\end{equation}
The element $u_\gamma$ should be regarded as an intertwiner
between representations with $\mc O (X)$-characters $x$ and
$\alpha_\gamma (x)$. We are interested in the finite type algebra
\begin{equation}\label{eq:1.13}
A\al = \mc O (X ; M_N (\mh C ))^\Gamma
\end{equation}
of $\Gamma$-invariant elements. We note that restriction of a
module from $A\al$ to $\mc O (X)^\Gamma$ defines a continuous
finite to one surjection \cite[Lemma 1]{KNS}
\begin{equation}\label{eq:1.3}
\theta : \mr{Prim}(A\al ) \to X / \Gamma \,.
\end{equation}
\textbf{Examples.}\\
Classical algebras of this type are
\begin{equation}\label{eq:1.14}
\begin{array}{lll}
\mc O (X)^\Gamma & = & \mc O (X / \Gamma ) \,, \\
\mc O \big( X ; \mr{End} (\mh C [\Gamma]) \big)^\Gamma
& \cong & \mc O (X) \rtimes \Gamma \,.
\end{array}
\end{equation}
For example, take $X = \mh C^3$ and $\Gamma = \mh Z / 3 \mh Z$,
acting through cyclic permutations of the coordinates. Put $A\al =
\mc O (X) \rtimes \Gamma$.  Almost all points $\Gamma x \in X /
\Gamma$ correspond to a unique irreducible $A\al$-module,
namely $\mr{Ind}_{\mc O (X)}^{A\al} \mh C_x$. Only the points
$(z,z,z)$ with $z \in \mh C$ carry three irreducible $\mc O (X)
\rtimes \Gamma$-modules, of the form $\mh C_{(z,z,z)} \otimes \mh
C_{\zeta}$ with $\zeta$ a cubic root of unity.
\\[2mm]
More generally, suppose that we have a larger group $G$ with a
normal subgroup $N$ such that $\Gamma = G / N$. Let $(\pi ,V)$ be a
$G$-representation on which $N$ acts by a character. Then
\[
(g \cdot f)(x) = \pi (g) f (\alpha_{g N}^{-1} x) \pi (g^{-1})
\]
defines an action of $G$ on $\mc O (X ; \mr{End} (V))$ which
factors through $\Gamma$, so
\[
\mc O (X ; \mr{End} (V))^G = \mc O (X ; \mr{End} (V))^\Gamma \,.
\]
If we put $u_\gamma = \pi (g)$ for some $g$ with $g N = \gamma$
then we are in the setting of \eqref{eq:1.2}. Yet in general
there is no canonical choice for $u_\gamma$, and we end up with a
nontrivial cocycle $\lambda$. (In fact this a typical example of a
projective $\Gamma$-representation.)\\[2mm]

The natural Fr\'echet completion of \eqref{eq:1.13} is
\begin{equation}\label{eq:1.22}
A\sm = C^\infty (X ; M_N (\mh C))^\Gamma \,.
\end{equation}
This algebra has the same spectrum as $A\al$, but the two algebras 
induce different Jacobson topologies on this set. The Jacobson 
topology from $A\sm$ is finer, and makes Prim$ (A\sm )$ a non-Hausdorff manifold. 
(By this we mean a second countable topological space in which every point has a
neighborhood that is homeomorphic to $\mh R^n$.)

The map \eqref{eq:1.4} is an isomorphism in the special cases where 
$A\al$ is as in \eqref{eq:1.14} and $A\sm$
as in \eqref{eq:1.22}, as follows from comparing \cite{KNS} and
\cite{Was}. We will show that it holds much more generally, for
example if we take $A\sm = C^\infty (X' ; M_N (\mh C))^\Gamma$
with $X'$ a suitable deformation retract of $X$. Such an algebra
is finitely generated as a $C^\infty (X')^\Gamma$-module, and
therefore we will call it a (topological) finite type algebra.
\vspace{4mm}

\section{The commutative case}
\label{sec:1.2}

Let $CC_{**}(A)$ be the periodic cyclic bicomplex associated to an algebra $A$
\cite[Section 1]{Lod}. Its terms are of the form $A^{\otimes n}$ and its homology
is by definition $HP_* (A)$.
For topological algebras we must specify which particular topological
tensor product we wish to use in cyclic theory. By default we work with a
Fr\'echet algebra $A$ whose topology is defined by
submultiplicative seminorms, and with the completed projective
tensor product $\hot$. Since this is a completion of the algebraic
tensor product we get natural maps
\[
\begin{array}{lll}
CC_{**}(A , \otimes) & \to & CC_{**}(A ,\hot) \,, \\
HP_* (A , \otimes) &\to & HP_* (A ,\hot) \,.
\end{array}
\]
We abbreviate $CC_{**}(A) = CC_{**}(A ,\hot)$ and $HP_* (A) = HP_* (A ,\hot )$ for Fr\'echet
algebras. Recall that an extension of topological algebras is admissible if it is split 
exact in the category of topological vector spaces. An ideal $I$ of $A$ is admissible if\\ 
$0 \to I \to A \to A/I \to 0$ is admissible. 

Any extension $0 \to A \to B \to C \to 0$ gives rise to a short exact sequence of
differential complexes
\begin{align*}
& 0 \to CC_{**}(B,A) \to CC_{**}(B) \to CC_{**}(C) \to 0 \,, \\
& CC_{**}(B,A) := \ker \big( CC_{**}(B) \to CC_{**}(C) \big) \,.
\end{align*}
By a standard construction in homological algebra this leads to a long exact sequence
\begin{equation}\label{eq:1.19}
\cdots \to HP_i (B,A) \to HP_i (B) \to HP_i (C) \to HP_{i+1}(B,A) \to \cdots
\end{equation}
Actually this sequence wraps up to an exact hexagon, because $HP_{i+2} \cong HP_i$.
The inclusion $CC_{**}(A) \to CC_{**}(B,A)$ induces a map $HP_* (A) \to HP_* (B,A)$.
One of our main tools will be the excision property of \pch \!:

\begin{thm}\label{thm:1.2}
Let $0 \to A \to B \to C \to 0$ be an extension of \textup{(}nontopological\textup{)}
algebras or an admissible extensions of Fr\'echet algebras. Then $HP_* (A) \to HP_* (B,A)$
is an isomorphism, and \eqref{eq:1.19} yields an exact hexagon
\[
\hexagon{HP_0 (A)}{HP_0 (B)}{HP_0 (C)}{HP_1 (C)}{HP_1 (B)}{HP_1 (A)}
\]
\end{thm}
\emph{Proof.}
The basic version of this theorem is due to Wodzicki \cite{Wod}.
It was proved in general by Cuntz and Quillen \cite{CuQu,Cun}.
$\qquad \Box$ \\[2mm]

Suppose that we want to prove that an algebra homomorphism
$\phi : A \to B$ induces an isomorphism on $HP_*$. The excision
property can be used as follows:

\begin{lem}\label{lem:1.3}
Suppose that there are sequences of ideals
\begin{align*}
& A = \, I_0 \supset \, I_1 \supset \cdots \supset I_d = 0 \\
& B = J_0 \supset J_1 \supset \cdots \supset J_d = 0
\end{align*}
with the properties
\begin{itemize}
\item $\phi (I_p ) \subset J_p$ for all $p \geq 0$,
\item $HP_* (I_{p-1} / I_p ) \to HP_* (J_{p-1} / J_p)$ is an
isomorphism for all $p \geq 0$,
\item if $B$ \textup{(}respectively $A$\textup{)} is Fr\'echet then the ideals
$J_p$ \textup{(}respectively $I_p$\textup{)} are admissible.
\end{itemize}
Then $HP_* (\phi) : HP_* (A) \to HP_* (B)$ is an isomorphism.
\end{lem}
\emph{Proof.}
Left as an exercise. Use the five lemma. $\qquad \Box$
\\[2mm]

Generally speaking a good tool to compute the \pch of a
finite type algebra $A$ is a filtration by ideals $I_p$ such that
the successive quotients $I_{p-1} / I_p$ behave like
commutative algebras. In particular Prim$ (I_{p-1}/I_p )$
should be a (separated) affine variety, so this gives rise a
kind of stratification of Prim$ (A)$. This can be formalized with the notion of 
an abelian filtration \cite{KNS}.

To describe suitable smooth analogues of $A\al$ we must say what
precisely we mean by smooth functions on spaces that are not
manifolds. Let $Z \subset Y$ be subsets of a smooth manifold $X$ and let
$V$ be a complete topological vector space.
\[
\begin{array}{lll}
C^\infty (Y ; V) & := & \big\{ f : Y \to V \,|\, \exists \text{
open } U \subset X , \tilde f \in C^\infty (U ; V) :
Y \subset U , \tilde f \big|_Y = f \big\} \\
C^\infty_0 (Y,Z) & := &
\{ f \in C^\infty (Y ; \mh C) : f \big|_Z = 0 \} \\
C^\infty_0 (Y,Z;V) & := & \{ f \in C^\infty (Y ; V) : f \big|_Z = 0 \}
\end{array}
\]
Recall that a corner in a manifold is a point that has a neighborhood 
homeomorphic to $\mh R^n \times [0,\infty )^m$, with $m > 0$. To 
apply excision we will often need the following result of Tougeron:

\begin{thm}\label{thm:1.9} \textup{\cite[Th\'eor\`eme IX.4.3]{Tou}}\\
Let $Y$ be a smooth manifold and $Z$ a smooth submanifold, both
possibly with corners. The following extension is admissible:
\[
0 \to C_0^\infty (Y,Z) \to C^\infty (Y) \to C^\infty (Z) \to 0 \,.
\]
\end{thm}

For completeness we include an extended version of the Hochschild--Kostant--Rosenberg
theorem for \pch \!. We abbreviate
\[
H^{[n]}(Y) = {\ts \bigoplus_{m \in \mh Z}} \check H^{n+2m}(Y ; \mh C) \,.
\]
\begin{thm}\label{thm:1.8}
Let $Y$ be a smooth manifold, possibly noncompact and with corners.
There is a natural isomorphism 
\[
HP_* (C^\infty (Y)) \cong H^{[*]}(Y ) \,.
\]
\end{thm}
\emph{Proof.} 
For $Y$ compact and without boundary this is due to Connes \cite[p. 130]{Con}, who in fact
proved the much stronger statement
\begin{equation}\label{eq:1.15}
HH_* (C^\infty (Y)) \cong \Omega^* (Y) \,.
\end{equation}
Here $HH_*$ denotes Hochschild homology and $\Omega^*$ means differential forms with 
complex values. By Corollary 4.3 and Theorem 7.1 of \cite{BLT} \eqref{eq:1.15} still holds if
$Y$ is allowed to have corners and may be noncompact. Hence $HP_* (C^\infty (?))$ and 
$H^{[*]} (? )$ agree at least locally. By Theorems \ref{thm:1.2} and \ref{thm:1.9} both
these functors satisfy excision, so they agree on every manifold. $\qquad \Box$
\\[2mm]

Let $Y$ be a complex affine variety and $Z$ a closed subvariety,
both possibly reducible and singular. In line with the above we write
\[
\begin{array}{lll}
\mc O_0 (Y,Z) & := & \{ f \in \mc O (Y) : f \big|_Z = 0 \} \,, \\
\mc O_0 (Y,Z;V) & := & O_0 (Y,Z) \otimes V \,.
\end{array}
\]
Let $\check H^* (Y,Z;\mh C )$ denote the \v Cech cohomology of the pair
$(Y,Z)$, with complex coefficients and with respect to the analytic
topology. Because $HP_*$ and $K_*$ have a $\mh Z / 2 \mh Z$-grading,
it is convenient to impose this also on \v Cech cohomology. Therefore we write
\begin{equation}\label{eq:1.6}
H^{[n]} (Y,Z) := {\ts \bigoplus_{m \in \mh Z}} \check H^{n+2m} (Y,Z;\mh C )
\end{equation}
Now we are ready to state and prove the comparison theorem
for the \pch of commutative algebras.

\begin{thm}\label{thm:1.4}
\begin{description}
\item[a)] There is a natural isomorphism
$HP_* \big( \mc O_0 (Y,Z) \big) \cong H^{[*]} (Y,Z)$.
\item[b)] Suppose that $Y \setminus Z$ is nonsingular and that
$\tilde C_0^\infty (Y,Z)$ is a Fr\'echet algebra with the properties
\item[$\bullet$] $\mc O_0 (Y,Z) \subset \tilde C_0^\infty (Y,Z) \subset
C_0^\infty (Y,Z)$ \,,
\item[$\bullet$] if the partial derivatives of $f \in C_0^\infty (Y,Z)$ all
vanish on $Z$ then $f \in \tilde C_0^\infty (Y,Z)$.\\
Then $HP_* (\mc O_0 (Y,Z)) \to HP_* \big( \tilde C_0^\infty (Y,Z) \big)$ 
is an isomorphism.
\end{description}
\end{thm}
\emph{Proof.}
a) was proved in \cite[Theorem 9]{KNS}.\\
b) By assumption $Y\setminus Z$ with the analytic topology is a smooth manifold. 
Let $N$ be a closed neighborhood of $Z$ in $Y$, such that $Y \setminus N$ is a 
smooth manifold and a deformation retract of $Y \setminus Z$. Let 
$r : [0,1] \times Y \to Y$ be a smooth map with the properties
\begin{enumerate}
\item $r_1 = \mr{id}_Y$ where $r_t (y) := r(t,y) $,
\item $r_t (z) = z \; \forall z \in Z , t \in [0,1] $,
\item $r_t (N) \subset N \; \forall t \in [0,1]$ and $r_0 (N) = Z $,
\item $r_t^{-1} (Z)$ is a neighborhood of $Z \; \forall t < 1 $,
\item $r_t \colon r_t^{-1}(Y \setminus Z) \to Y \setminus Z$ is a
diffeomorphism $\forall t \in [0,1] $.
\end{enumerate}
Consider the algebra homomorphisms
\[
\begin{array}{lll@{\qquad}lll}
C_0^\infty (Y,N) & \to & \tilde C_0^\infty (Y,Z), & f & \to & f \,, \\
\tilde C_0^\infty (Y,Z) & \to & C_0^\infty (Y,N), &
f & \to & f \circ r_0 \,.
\end{array}
\]
By construction these are diffeotopy equivalences. The diffeotopy
invariance of $HP_* (?,\hot )$ \cite[p. 125]{Con} ensures that
\[
HP_* \big( \tilde C_0^\infty (Y,Z) \big) \cong
HP_* \big( C_0^\infty (Y,N) \big) \,.
\]
Write $\tilde Y = r_{1/2}^{-1} (Y \setminus Z)$ and
$\tilde N = N \cap \tilde Y$. This could look like
\\[1mm]
\includegraphics[width=8cm,height=5cm]{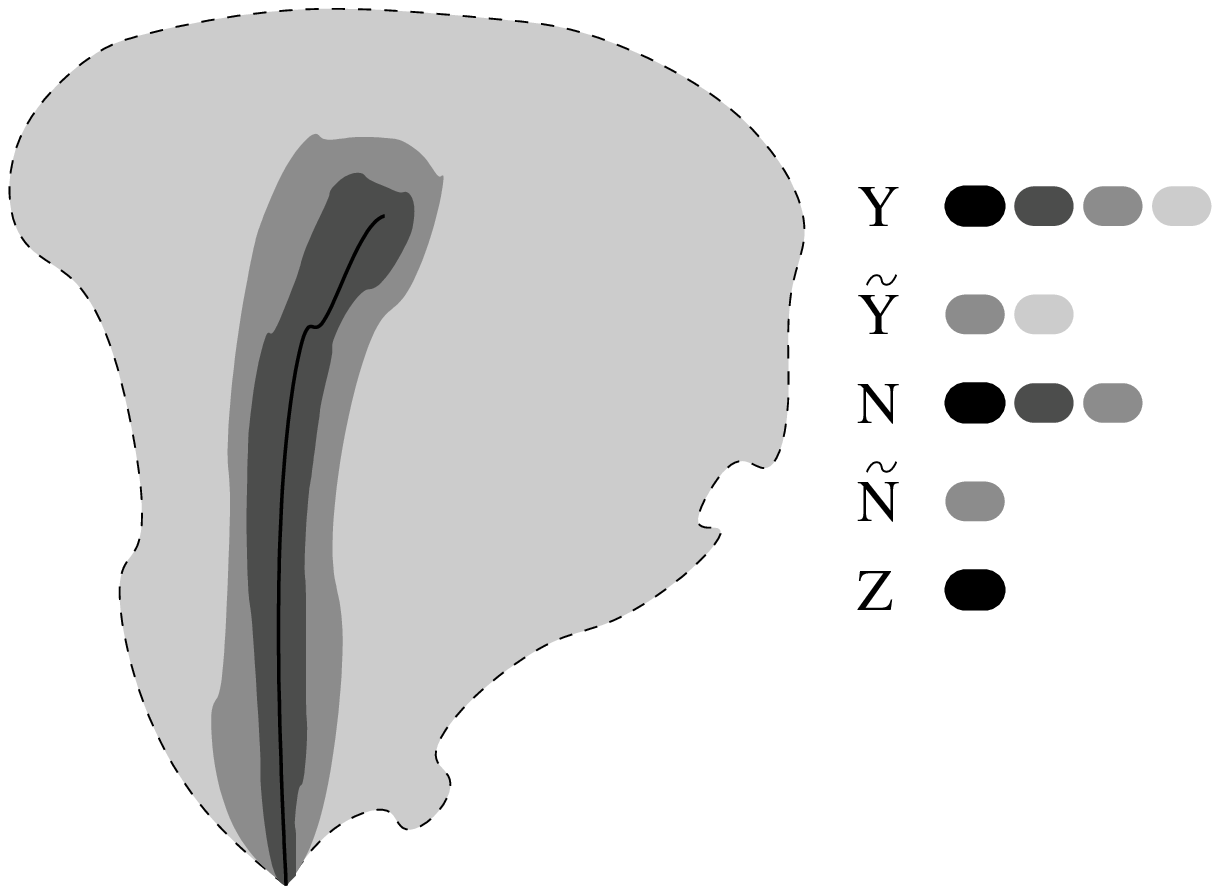}

\noindent By Theorem \ref{thm:1.9} there is an admissible extension
\[
0 \to C_0^\infty (Y,N) = C_0^\infty \big( \tilde Y , \tilde N \big)
\to C^\infty \big( \tilde Y \big) \to C^\infty \big( \tilde N \big)
\to 0 \,.
\]
Combining this with Theorems \ref{thm:1.2} and \ref{thm:1.8} yields natural
isomorphisms
\[
HP_* \big( C^\infty (\tilde Y ,\tilde N) \big) \cong H^{[*]} \big(
\tilde Y ,\tilde N \big) \cong H^{[*]} (Y,N) \cong H^{[*]} (Y,Z) \,.
\]
Now consider the diagram
\[
\begin{array}{ccc}
HP_* (\mc O_0 (Y,Z)) & \cong & H^{[*]} (Y,Z) \\
\downarrow & & || \\
HP_* \big( C^\infty (\tilde Y ,\tilde N) \big) & \cong &
H^{[*]} (Y,Z)
\end{array}
\]
It commutes by naturality, so the arrow is an isomorphism.
$\qquad \Box$
\vspace{4mm}

\section{Comparison with topological $K$-theory}
\label{sec:1.3}

Let $\Gamma$ be a finite group which acts by diffeomorphisms on a 
smooth manifold $X$. We will frequently meet algebras of the form
\begin{equation}\label{eq:1.5}
C_0^\infty (Y,Z;M_N (\mh C) )^\Gamma
\end{equation}
where $Y$ and $Z$ are $\Gamma$-stable closed submanifolds of $X$.
We allow our manifolds to have corners (and in particular a
boundary), since these appear naturally in orbifolds. But we must
be careful, because the algebra $C^\infty (X)^\Gamma$ of smooth
functions on the orbifold $X/ \Gamma$ does not contain all smooth
functions on the manifold $X / \Gamma$. Namely, there are some
conditions for the partial derivatives of $f \in C^\infty
(X)^\Gamma$ at the corners.

This makes it rather tricky to compute the \pch of algebras
like \eqref{eq:1.5}. Actually such algebras would look much
better if we could replace smooth functions by continuous
functions, because continuous functions are not bothered by
mild singularities like corners. But then another problem
pops up, that $HP_*$ tends to give tautological results
for Banach algebras. For example, if $K$ is any compact
Hausdorff space then
\begin{align*}
& HP_0 (C(K)) = C(K)\,, \\
& HP_1 (C(K)) = 0 \,.
\end{align*}
To overcome these inconveniences we will compute the \pch of \eqref{eq:1.4} via 
its topological $K$-theory. The $K$-theory of Fr\'echet algebras was
defined by Phillips \cite{Phi}. As is well-known, these theories are related 
by a Chern character $ch : K_* \to HP_*$.

\begin{thm}\label{thm:1.5} \textup{\cite[Theorem 16]{Nis}}\\
Let $0 \to A \to B \to C \to 0$ be an extension of Fr\'echet
algebras. The various Chern characters form a commutative
diagram
\[
\begin{array}{ccccccccccc}
K_1(A) & \to & K_1(B) & \to & K_1(C) & \to & K_0(A) & \to &
 K_0(B) & \to & K_0(C) \\
\downarrow & & \downarrow & & \downarrow & & \downarrow & &
\downarrow & & \downarrow \\
HP_1(A) & \to & HP_1(B) & \to & HP_1(C) & \to & HP_0(A) & \to &
 HP_0(B) & \to & HP_0(C)
\end{array}
\]
If the extension is admissible and $\eta : K_0 (C) \to K_1 (A)$
and $\partial : HP_0 (C) \to HP_1 (A)$ denote the connecting
maps, then $ch \circ \eta = 2\pi i \, \partial \circ ch$.
\end{thm}
\vspace{2mm}

Let $\mc{CIA}$ be the class of Fr\'echet algebras $A$ for which
\begin{equation}
ch \otimes \mr{id}_{\mh C} : K_* (A) \otimes_{\mh Z} \mh C \to HP_*
(A) \otimes_{\mh C} \mh C = HP_* (A)
\end{equation}
is an isomorphism.

\begin{cor}\label{cor:1.10}
Let $0 \to A \to B \to C \to 0$ be an admissible extension of
Fr\'echet algebras. If two of ${A,B,C}$ belong to the class
$\mc{CIA}$, then so does the third.
\end{cor}
\emph{Proof.}
This follows from Theorems \ref{thm:1.2} and \ref{thm:1.5},
in combination with Bott periodicity and the five lemma.
$\qquad \Box$ \\[3mm]

The class $\mc{CIA}$ is very large, since it is also closed
under countable direct products, tensoring with $M_n (\mh C)$
and diffeotopy equivalences. Furthermore all topological
finite type algebras are in $\mc{CIA}$:

\begin{thm}\label{thm:1.6}
Let $\Gamma$ be a finite group acting \textup{(}by $\alpha$\textup{)} 
on a smooth manifold $X$ and let $u_\gamma \in GL_N (C^\infty (X))$
be elements such that
\[
\gamma \cdot f = u_\gamma (f \circ \alpha_\gamma^{-1}) u_\gamma^{-1}
\]
defines an action of $\Gamma$ on $C^\infty (X ; M_N (\mh C))$.
Let $Z \subset Y$ be $\Gamma$-stable submanifolds of $X$, possibly
with corners. Then $C_0^\infty (Y,Z;M_N (\mh C))^\Gamma$ belongs
to the class $\mc{CIA}$.
\end{thm}
\emph{Proof.}
According to Theorem \ref{thm:1.9}
\begin{equation}\label{eq:1.7}
0 \to C_0^\infty (X,Y) \to C^\infty (X) \to C^\infty (Y) \to 0
\end{equation}
is an admissible extension. Hence so is
\begin{equation}\label{eq:1.8}
0 \to C_0^\infty (X,Y;M_N (\mh C )) \to C_0^\infty (X,Z;M_N (\mh C
)) \to C_0^\infty (Y,Z;M_N (\mh C )) \to 0 \,.
\end{equation}
Because $\Gamma$ is finite the same holds for the subalgebras of
$\Gamma$-invariants in \eqref{eq:1.7} and \eqref{eq:1.8}. Together
with Corollary \ref{cor:1.10} this reduces the proof to the case
$Z = \es$.

Thus we have to show that $C^\infty (Y;M_N (\mh C ))^\Gamma$
is in $\mc{CIA}$. Except for a detail this is the content of
\cite[Theorem 6]{Sol1}. The small complication is that in \cite{Sol1}
the author considered only $\Gamma$-manifolds $Y$ without corners,
because according to \cite{Ill} those have smooth equivariant
triangulations. However, $Y$ admits such a triangulation even if it has
corners, because it is embedded in the smooth $\Gamma$-manifold $X.
\qquad \Box$
\vspace{4mm}

\section{The general case}
\label{sec:1.4}

First we discuss a motivating concept for our comparison theorem.
Suppose that $\Gamma$ acts on $\mh Z^n$, and consider the tori
\[
\begin{array}{lllllll}
X & := & \mr{Hom}_{\mh Z}(\mh Z^n ,\mh C^\times) & \cong &
\big( \mh C^\times \big)^n & \cong &
\mr{Prim} \, \big( \mh C [\mh Z^n ] \big) \,, \\
X' & := & \mr{Hom}_{\mh Z}(\mh Z^n ,S^1) & \cong & \big( S^1 \big)^n
& \cong &\mr{Prim} \, \big( \mc S (\mh Z^n )\big) \,,
\end{array}
\]
where the $\mc S$ stands for complex valued Schwartz functions. We
want to compare the \pch of the algebras
\[
\begin{array}{lcccr}
A\al & = & \mh C [\mh Z^n ] \rtimes \Gamma & = & \mc O (X) \rtimes
\Gamma \,, \\
A\sm & = & S (\mh Z^n ) \rtimes \Gamma & = & C^\infty (X') \rtimes
\Gamma \,.
\end{array}
\]
Although these algebras definitely have different spectra, it is
natural to expect that $HP_* (A\al ) \cong HP_* (A\sm )$.
The best notion to explain this appears to be ``diffeotopy equivalence of
non-Hausdorff spaces". This is a typically noncommutative geometric
concept that might contradict one's intuition. The idea is that Prim$ (A\al )$
and Prim$ (A\sm )$ are equivalent in this specific sense, and for that very 
reason these algebras have the same \pch \!. 

Since usual homotopies do not see non-Hausdorff phenomena we have
to be careful in defining this notion. We say that a continuous map
$X \to Y$ is a homotopy (diffeotopy) equivalence of non-Hausdorff
spaces if there exist finite length stratifications of $X$ and $Y$ such that:
\begin{itemize}
\item all the strata are Hausdorff spaces,
\item the maps are compatible with the stratifications,
\item the induced maps on the strata are homotopy (diffeotopy) equivalences.
\end{itemize}
Notice the we do not require the existence of a continuous map from $Y$ to 
$X$, because that would exclude many interesting cases.
For example consider the plane with a doubled origin. It is contractible in the usual
sense, but as a non-Hausdorff space it is diffeotopy equivalent to two points!

Generally speaking an algebra homomorphism that induces a diffeotopy 
equivalence on primitive ideal spectra (endowed with a suitable ``analytic"
topology) should yield an isomorphism on periodic cyclic homology. Our notion
is probably not strong enough to prove things with, but it does provide
a generalization of Theorem \ref{thm:1.1} at a conceptual level. 

Thus inspired we require the following conditions for our
comparison theorem. Let $\Gamma$ be a finite group acting (by
$\alpha$) on a nonsingular complex affine variety $X$. Suppose that
we have elements $u_\gamma \in GL_N (\mc O (X))$ such that
\[
\gamma \cdot f = u_\gamma (f \circ \alpha_\gamma^{-1} ) u_\gamma^{-1}
\]
defines an action of $\Gamma$ on the algebra $\mc O (X;M_N (\mh C ))$.
Let $X'$ be a submanifold of $X$ with the following properties:
\begin{itemize}
\item $X'$ is smooth, but may have corners,
\item $X'$ is stable under the action of $\Gamma$,
\item the inclusion $X' \to X$ is a diffeotopy equivalence
in the category of smooth $\Gamma$-manifolds.
\end{itemize}
We write
\[
\begin{array}{lcr}
A\al & = & \mc O (X;M_N (\mh C ))^\Gamma \,, \\
A\sm & = & C^\infty (X';M_N (\mh C ))^\Gamma \,.
\end{array}
\]
The inclusion of Prim$ (A\sm )$ in Prim$ (A\al )$ is the prototype
of a diffeotopy equivalence of non-Hausdorff spaces.
\begin{thm}\label{thm:1.7}
The natural map $A\al \to A\sm$ induces an isomorphism
\[
HP_* (A\al ) \to HP_* (A\sm ) \,.
\]
\end{thm}
\emph{Proof.} We will use Lemma \ref{lem:1.3} to reduce the proof to
manageable pieces. For every subset $H \subset \Gamma$ the variety
$X^H$ is nonsingular and $X'^H = X^H \cap X'$ is a submanifold. Let
$\mc L$ be the collection of all the irreducible components of all the $X^H$, with 
$H$ running over all subsets of $\Gamma$. Let $\mc L_p$ be its subset of 
elements of dimension $\leq p$ and define $\Gamma$-stable closed subvarieties
\[
X_p := {\ts \bigcup_{V \in \mc L_p}} \, V \,.
\]
By the third condition above $X'_p := X_p \cap X'$ is $\Gamma$-equivariantly 
diffeotopy equivalent to $X_p$.

By construction the singularities of $X_p$
are all contained in $X_{p-1}$. Moreover, because the action of $\Gamma$
is locally linearizable, these singularities are all normal crossings.
Hence we have for arbitrary subsets $G,H \subset \Gamma $:
\begin{equation}\label{eq:1.9}
\begin{array}{lll}
X^G \cap X^H & = & X^{G \cup H} \,, \\
\mc O_0 (X^G \cup X^H , X^G \cap X^H ) & = &
\mc O_0 (X^G , X^{G \cup H}) \oplus \mc O_0 (X^H , X^{G \cup H}) \,, \\
X'^G \cap X'^H & = & X'^{G \cup H} \,, \\
C^\infty_0 (X'^G \cup X'^H , X'^G \cap X'^H ) & = &
C^\infty_0 (X'^G , X'^{G \cup H}) \oplus C^\infty_0 (X'^H , X'^{G \cup H}) \,.
\end{array}
\end{equation}
Consider the sequences of ideals
\begin{equation}\label{eq:1.10}
\begin{array}{lllll}
A\al & = & \, I_0 \supset \, I_1 \supset \cdots \supset I_{\dim X} & = & 0 \,, \\
A\sm & = & J_0 \supset J_1 \supset \cdots \supset J_{\dim X} & = & 0 \,, \\
I_p & = & \{ a \in A\al : a \big|_{X_p} = 0 \} & = &
\mc O_0 (X,X_p ;M_N (\mh C ))^\Gamma \,, \\
J_p & = & \{ a \in A\sm : a \big|_{X'_p} = 0 \} & = &
C^\infty_0 (X',X'_p ;M_N (\mh C ))^\Gamma \,.
\end{array}
\end{equation}
We want to compare the \pch of the quotients
\[
\begin{array}{llr}
I_{p-1}/I_p & \cong & \mc O_0 (X_p ,X_{p-1} ;M_N (\mh C ))^\Gamma \,, \\
J_{p-1}/J_p & \cong & C^\infty_0 (X'_p ,X'_{p-1} ;M_N (\mh C ))^\Gamma \,.
\end{array}
\]
Let $Z (B)$ denote the center of an algebra $B$. To the filtration \eqref{eq:1.10}
we associate the spaces
\begin{equation}\label{eq:1.20}
\begin{aligned}
& Y_p = \mr{Prim}\big(Z (A\al / I_p )\big) \,,\\
& Z_p = \{ I \in Y_p : Z (I_{p-1}/I_p) \subset I \} \,.
\end{aligned}
\end{equation}
The $Y_p$ are called the centers of the filtration, and the $Z_p$ the subcenters.
Notice that, unlike Prim$ (I_p )$, these are separated algebraic varieties. 
By \cite[Theorem 9]{KNS} there are natural isomorphisms 
\begin{equation}
HP_* \big( Z (I_{p-1}/I_p )\big) \cong H^{[*]} (Y_p ,Z_p ) \cong
HP_* \big( \mc O_0 (Y_p ,Z_p ) \big) \,.
\end{equation}
We claim that 
\begin{equation}\label{eq:1.16}
Z (I_p / I_{p-1} ) \to I_p / I_{p-1} 
\end{equation}
is a spectrum preserving morphism of finite type $\mc O (X)$-algebras.
To see this, we first consider the composite map
\begin{equation}\label{eq:1.17} 
\theta_p : \mr{Prim}(I_{p-1} / I_p ) \to \mr{Prim}(Z (I_p / I_{p-1} )) = 
Y_p \setminus Z_p \to (X_p \setminus X_{p-1} ) / \Gamma \,.
\end{equation}
For $x \in X_p \setminus X_{p-1}$ the image of 
\[
I_{p-1} / I_p \to  M_N (\mh C ) : f \mapsto f(x)
\]
is the semisimple algebra $S_x = \mr{End}_{\pi_x (\Gamma_x )} ( \mh C^N )$, where 
$\pi_x (\gamma ) = u_\gamma (x)$. Since $u_\gamma \in GL_N (\mc O (X))$ the type of 
$(\pi_x ,\mh C^N )$ as a projective $\Gamma_x$-representation cannot change along the 
irreducible components of $X^{\Gamma_x}$. Together with \eqref{eq:1.10} this implies 
that locally on $(X_p \setminus X_{p-1} ) / \Gamma \,, I_{p-1} / I_p$ is of the form 
$S_x \otimes \mc O_0 (U,U')$. Hence $Z(I_{p-1} / I_p )$ is locally of the form 
$Z (S_x) \otimes \mc O_0 (U,U')$, which proves our claim about \eqref{eq:1.16}. 
Now we may apply Theorem \ref{thm:1.1}, which tells us that
\begin{equation}
HP_* \big( Z (I_{p-1}/I_p )\big) \to HP_* (I_{p-1}/I_p )
\end{equation}
is an isomorphism.

According to a very general extension theorem for smooth functions
\cite[Theorem 0.2.1]{BiSc} the ideals $J_p$ are admissible in
$A\sm$. Alternatively, this can be derived from
Theorem \ref{thm:1.9}, using \eqref{eq:1.9}. From \eqref{eq:1.9} we
also see that $J_{p-1}/J_p$ is a finite direct sum of algebras
of the form considered in Theorem \ref{thm:1.6}, so
\begin{equation}
ch \otimes \mr{id} : K_* (J_{p-1}/J_p ) \otimes_{\mh Z} \mh C
\to HP_* (J_{p-1}/J_p )
\end{equation}
is an isomorphism. Moreover $J_{p-1}/J_p$ is dense and
holomorphically closed in
\[
A_p := C_0 (X_p , X_{p-1};M_N (\mh C ))^\Gamma \,.
\]
The density theorem for $K$-theory \cite[Th\'eor\`eme A.2.1]{Bost}
tells us that the inclusion $J_{p-1}/J_p \to A_p$ induces an
isomorphism
\begin{equation}\label{eq:1.11}
K_* (J_{p-1}/J_p ) \to K_* ( A_p ) \,.
\end{equation}
The same arguments apply to the center of $J_{p-1}/J_p$ so there
are natural isomorphisms
\begin{equation}\label{eq:1.12}
\begin{array}{lll}
K_* \big( Z (J_{p-1}/J_p ) \big) \otimes_{\mh Z} \mh C & \to &
HP_* \big( Z (J_{p-1}/J_p ) \big) \,, \\
K_* \big( Z (J_{p-1}/J_p ) \big) & \to & K_* \big( Z(A_p ) \big) \,.
\end{array}
\end{equation}
The spectrum of the algebras in \eqref{eq:1.11} and
\eqref{eq:1.12} is $Y'_p \setminus Z'_p$ where
\[
\begin{array}{lllll}
Y'_p & = & \mr{Prim}\big( Z (A\sm / J_p ) \big) &
= & Y_p \cap \theta_p^{-1} (X'_p / \Gamma ) \,, \\
Z'_p & = & \{ J \in Y'_p : Z (J_{p-1} / J_p ) \subset J \} &
= & Z_p \cap \theta_p^{-1} (X'_p / \Gamma ) \,,
\end{array}
\]
with $\theta_p$ as in \eqref{eq:1.17}. Since the cardinality of $\theta_p^{-1}(x)$ 
is locally constant for $x \in X_p \setminus X_{p-1}$, any $\Gamma$-equivariant 
diffeotopy implementing the diffeotopy equivalence $X'_p \to X_p$ naturally gives 
rise to a diffeotopy for the inclusion map $(Y'_p ,Z'_p ) \to (Y_p ,Z_p )$. Therefore
\begin{equation}
\check H^n (Y_p ,Z_p ;\mh C ) \to \check H^n (Y'_p ,Z'_p ;\mh C )
\end{equation}
is an isomorphism for all $n$.

Furthermore $A_p$ is a finite direct sum of algebras of the form
$C_0 (Y,Z;M_k (\mh C) )$ with $Y$ a connected manifold.
The center of such an algebra is $C_0 (Y,Z)$, which clearly
is Morita-equivalent to the algebra itself. Hence
\[
Z (A_p ) = C_0 (Y'_p , Z'_p )
\]
and the inclusion map induces an isomorphism
\[
K_* \big( C_0 (Y'_p , Z'_p ) \big) \to K_* (A_p ) \,.
\]
Returning to the smooth level we note that
\[
Z (J_{p-1}/J_p ) := \tilde C_0^\infty (Y'_p ,Z'_p )
\]
satisfies the conditions 1 and 2 of Theorem \ref{thm:1.4}.
The proof of Theorem \ref{thm:1.4} yields a natural isomorphism
\begin{equation}\label{eq:1.21}
HP_* \big( \tilde C_0^\infty (Y'_p ,Z'_p ) \big) \cong
H^{[*]} (Y'_p , Z'_p ) \,.
\end{equation}
Combining all the above we get a diagram
\[
\begin{array}{ccccccc}
HP_* (I_{p-1}/I_p ) & \xleftarrow{(4)} & HP_* \big( Z (I_{p-1}/I_p )
\big) & \cong & HP_* \big( \mc O_0 (Y_p ,Z_p ) \big) & \cong &
H^{[*]} (Y_p , Z_p ) \\
\downarrow {\scs (8)} & & \downarrow {\scs (7)} & &
\downarrow {\scs (6)} & & \downarrow {\scs (5)} \\
HP_* (J_{p-1}/J_p ) & \xleftarrow{(3)} & HP_* \big( Z (J_{p-1}/J_p ) \big) &
\cong & HP_* \big( \tilde C_0^\infty (Y'_p , Z'_p ) \big) & \cong &
H^{[*]} (Y'_p ,Z'_p ) \\
\uparrow & & \uparrow & & \uparrow & & \uparrow \\
K_* (J_{p-1}/J_p ) & \xleftarrow{(2)} & K_* \big( Z (J_{p-1}/J_p ) \big) &
\cong & K_* \big( \tilde C_0^\infty (Y'_p , Z'_p ) \big) & \cong &
K^* (Y'_p ,Z'_p ) \\
\downarrow & & \downarrow & & \downarrow & & || \\
K_* (A_p ) & \xleftarrow{(1)} & K_* \big( Z (A_p ) \big) &
\cong & K_* \big( C_0 (Y'_p , Z'_p ) \big) & \cong & K^* (Y'_p ,Z'_p )
\end{array}
\]
that is commutative because all the maps are natural. So far we know that:
\begin{itemize}
\item the maps from row 3 to row 4 are isomorphisms by the density
theorem in topological $K$-theory,
\item the Chern characters from row 3 to row 2 become isomorphisms
after tensoring with $\mh C$,
\item (1), (4) and (5) are isomorphisms.
\end{itemize}
With some obvious diagram chases we first deduce that (2) and (3)
are isomorphisms, and then that (6), (7)
and finally (8) are isomorphisms. $\qquad \Box$
\\[2mm]

\noindent\textbf{Example.}\\
Let $X = \mh C^2 \,, X' = [-1,1]^2 \subset \mh R^2 \subset \mh
C^2$ and $\Gamma = \{ \pm 1 \}^2$. We describe the stratifications
of the spectra of the algebras
\[
\begin{array}{lcr}
A\al & = & \mc O (X) \rtimes \Gamma \,, \\
A\sm & = & C^\infty (X') \rtimes \Gamma \,.
\end{array}
\]
First the strata of $X$ and $X'$ :
\[
\begin{array}{lll@{\qquad}lll}
X_0 & = & \{ (0,0) \} \,, & X'_0 & = & \{ (0,0) \} \,,\\
X_1 & = & \{ 0 \} \times \mh C \: \cup \: \mh C \times \{ 0 \} \,, &
X'_1 & = & \{ 0 \} \times [-1,1] \: \cup \: [-1,1] \times \{ 0 \} \,, \\
X_2 & = & \mh C^2 \,, & X'_2 & = & [-1,1]^2 \,.
\end{array}
\]
Let $\sigma$ and $\tau$ be the two irreducible representations of the
group $\{ \pm 1 \}$. We extend them to representations $\sigma_0 , \tau_0$
of $C^\infty ([-1,1]) \rtimes \{ \pm 1 \}$ with central character $0 \in [-1,1]$.
The centers of the filtration are:
\[
\begin{array}{lll}
Y_0 & = & \{ \sigma_0 ,\tau_0 \}^2 \,, \\
Y_1 & = & \{ \sigma_0 ,\tau_0 \} \times \mh C / \{ \pm 1 \} \: \cup \:
\mh C / \{ \pm 1 \} \times \{ \sigma_0 ,\tau_0 \} \big/ \sim \\
& \cong & \{ 0 \} \times \mh C \: \cup \: \mh C \times \{ 0 \} \,, \\
Y_2 & = & X / \Gamma \: \cong \: ( \mh C / \{ \pm 1 \} )^2 \,, \\
Y'_0 & = & \{ \sigma_0 ,\tau_0 \}^2 \,, \\
Y'_1 & = & \{ \sigma_0 ,\tau_0 \} \times [-1,1] \: \cup \:
[-1,1] \times \{ \sigma_0 ,\tau_0 \} \big/ \sim \\
& \cong & \{ 0 \} \times [-1,1] \: \cup \: [-1,1] \times \{ 0 \} \,, \\
Y'_2 & = & X' / \Gamma \: \cong \: [0,1]^2 \,.
\end{array}
\]
where the equivalence relation $\sim$ identifies all the points lying over $(0,0) \in \mh C^2$.
Next we write down the subcenters of the filtration:
\[
\begin{array}{lll@{\qquad}lll}
Z_0 & = & \es \,, & Z'_0 & = & \es \,,\\
Z_1 & = & \{ (0,0) \} \,, & Z'_1 & = & \{ (0,0) \} \,, \\
Z_2 & = & X_1 \,, & Z'_2 & = & X'_1 \,.
\end{array}
\]
Finally we mention the primitive ideal spectra of the subquotients of the filtrations:
\[
\begin{array}{lll@{\qquad}lll}
(Y_0 \setminus Z_0 ) / \Gamma & \cong & 4 \text{ points,} & 
(Y'_0 \setminus Z'_0 ) / \Gamma & \cong & 4 \text{ points} \,, \\
(Y_1 \setminus Z_1 ) / \Gamma & \cong & \{ 1,2 \} \times \mh C^\times / \{\pm 1\} \,, &
(Y'_1 \setminus Z'_1 ) / \Gamma & \cong & \{1,2 \} \times (0,1] \,, \\
(Y_2 \setminus Z_2 ) / \Gamma & \cong & \big( \mh C^\times / \{\pm 1\} \big)^2 \,, & 
(Y'_2 \setminus Z'_2 ) / \Gamma & \cong & (0,1]^2 \,.
\end{array}
\]

\chapter{Some representation theory of reductive $p$-adic groups}

\section{Convolution algebras}
\label{sec:2.1}

In this chapter we collect some important results concerning
smooth representations of reductive $p$-adic groups. Good sources
for the theory discussed here are
\cite{BeDe,Car,Sil2,SSZ,Tits,Wal}.

Let $\mh F$ be a non-Archimedean local field with discrete valuation $v$ 
and norm $\norm{.}_{\mh F}$. We assume that the cardinality of the residue 
field is a power $q$ of a prime $p$. Let $\mc G$ be a connected reductive algebraic 
group defined over $\mh F$, and let $G = \mc G (\mh F )$ be the group of 
$\mh F$-rational points. We briefly call $G$ a reductive $p$-adic group.

We denote the collection of compact open subgroups of $G$ by
CO$(G)$. A representation $V$ of $G$ is called smooth if every $v
\in V$ is fixed by a compact open subgroup, or equivalently if $V
= \cup_{K \in \mr{CO}(G)} V^K$. We say that such a smooth representation $V$ 
is admissible if every $V^K$ has finite dimension. For example every smooth 
$G$-representation of finite length is admissible \cite[3.12]{BeDe}.
Let Rep$ (G)$ be the category of
smooth $G$-representations on complex vector spaces, and let Irr$ (G)$ be 
the set of equivalence classes of irreducible objects in Rep$ (G)$.
The Jordan--H\"older content JH$ (V)$ is the collection of all elements of Irr$ (G)$
which are equivalent to a subquotient of the $G$-representation $V$.

Fix a Haar measure $d\mu$ on $G$. Recall that the convolution
product of two functions $f,f' : G \to \mh C$ is defined as
\[
(f * f')(g') = \int_G f(g) f'(g^{-1} g') \, \textup{d}\mu (g) \,.
\]
For $K \in \mr{CO}(G)$ we let $\mc H (G,K)$ be the convolution
algebra of $K$-biinvariant complex-valued compactly supported
functions on $G$. This is called the Hecke algebra of $(G,K)$. Our
main subject of study will be the Hecke algebra of $G$, which
consists of all compactly supported locally constant functions on $G$:
\begin{equation}
\mc H (G) := {\ts \bigcup_{K \in \mr{CO}(G)} } \mc H (G,K) \,.
\end{equation}
For every $K \in \mr{CO}(G)$ there is an idempotent $e_K \in \mc H
(G)$, which is $\mu (K)^{-1}$ times the characteristic function of
$K$. Notice that
\begin{equation}
\mc H (G,K) = e_K \mc H (G) e_K = \mc H (G)^{K \times K} ,
\end{equation}
where $G \times G$ acts on $\mc H (G)$ by left and right
translations. In particular, the nonunital algebra $\mc H (G)$ is
idempotented, which assures that many properties of unital
algebras also hold for $\mc H (G)$.

An $\mc H (G)$-module $V$ is called nondegenerate or essential if
\begin{equation}\label{eq:2.1}
\mc H \cdot V = V ,
\end{equation}
or equivalently if for all $v \in V$ there exists a $K \in \mr{CO}(G)$
such that $e_K \cdot v = v$. A smooth $G$-representation $(\pi ,V)$ is
made into an essential $\mc H (G)$-module by
\begin{equation}\label{eq:2.12}
\pi (f) v = \int_G f (g) \pi (g) v \, \textup{d}\mu (g) \qquad v \in V , f \in \mc H (G) \,.
\end{equation}
This leads to an equivalence between Rep$(G)$ and the category of
essential $\mc H (G)$-modules. Hence we may identify the 
primitive ideal spectrum of $\mc H (G)$ with Irr$ (G)$.

Let $\mc S (G,K)$ be the space of rapidly decreasing
$K$-biinvariant functions on $G$. According to \cite[Theorem
29]{Vig} this is a unital nuclear Fr\'echet *-algebra.
Harish-Chandra's Schwartz algebra consists of all uniformly
locally constant rapidly decreasing functions on $G$:
\begin{equation}
\mc S (G) := {\ts \bigcup_{K \in \mr{CO}(G)} } \mc S (G,K) \,.
\end{equation}
Endowed with the inductive limit topology this is a complete locally convex 
topological algebra with separately continuous multiplication. Clearly
\begin{equation}
\mc S (G,K) = e_K \mc S (G) e_K = \mc S (G)^{K \times K} .
\end{equation}
If $(K_i )_{i=1}^\infty $ is a decreasing sequence of compact open
subgroups of $G$ which forms a neighborhood basis of the unit
element $e \in G$, then
\begin{equation}\label{eq:2.2}
\mc S (G) = {\ts \bigcup_{i=1}^\infty } \mc S (G,K_i )
\end{equation}
is a strict inductive limit of nuclear Fr\'echet spaces. Nevertheless
$\mc S (G)$ is not metrizable. We say that a smooth
$G$-representation $(\pi, V)$ is tempered if the $\mc H (G)$-module structure
extends to $\mc S (G)$. If $(\pi,V)$ is admissible, then there is at most one
such extension, see \cite[p. 51]{SSZ}. Thus we have
\begin{itemize}
\item the category $\mr{Rep}^t (G)$ of tempered smooth
$G$-representations,
\item the space $\mr{Irr}^t (G)$ of equivalence classes of
irreducible objects in $\mr{Rep}^t (G) $,
\item the primitive ideal spectrum of $\mc S (G)$, which by
\cite[p. 52]{SSZ} can be identified with $\mr{Irr}^t (G)$.
\end{itemize}
Furthermore we consider the reduced $C^*$-algebra of $G$. By
definition $C_r^* (G)$ is the completion of $\mc H (G)$ with
respect to the operator norm coming from the left regular
representation of $G$ on $L^2 (G)$. For $K \in \mr{CO}(G)$ let
$C_r^* (G,K)$ be the norm closure of $\mc H (G,K)$ in $B(L^2 (G))$. 
This is a unital type I $C^*$-algebra which contains $\mc S
(G,K)$ as a holomorphically closed dense subalgebra 
\cite[Theorem 29]{Vig}. Moreover by \cite[p. 53]{SSZ}
\begin{equation}
C_r^* (G,K) = e_K C_r^* (G) e_K = C_r^* (G)^{K \times K} .
\end{equation}
Therefore we can construct $C_r^* (G)$ also as an inductive limit
of $C^*$-algebras:
\begin{equation}\label{eq:2.14}
C_r^* (G) = \varinjlim_{K \in \mr{CO}(G)} C_r^* (G,K) \,.
\end{equation}
Having introduced these algebras we will describe the Bernstein
decomposition of Rep$ (G)$.
Suppose that $P$ is a parabolic subgroup of $G$ and that $P = M 
\ltimes N$ where $N$ is the unipotent radical of $P$ and $M$ is a Levi
subgroup. Although $G$ and $M$ are unimodular the modular function
$\delta_P$ of $P$ is general not constant. To be precise
\begin{equation}
\delta_P (m n) = \norm{\det \big( \mr{ad}(m) \big|_{\mf n} \big)
}_{\mh F} \qquad m \in M , n \in N
\end{equation}
where $\mf n$ is the Lie algebra of $N$. For $\sigma \in \mr{Rep}(M)$ one defines
\[
I_P^G (\sigma ) := \mr{Ind}_P^G (\delta_P^{1/2} \otimes \sigma ) \,.
\]
This means that we first inflate $\sigma$ to $P$, then we twist it with $\delta_P^{1/2}$ 
and finally we take the smooth induction to $G$. The twist is useful to preserve unitarity. 
The functor $I_P^G$ is known as parabolic induction. It is transitive in the sense that 
for any parabolic subgroup $Q \subset P$  we have 
\[
I_Q^G = I_P^G \circ I^M_{Q \cap M} \,.
\]
Let $\sigma$ be an irreducible supercuspidal representation of $M$. Thus
the restriction of $\sigma$ to the derived group of $M$ is unitary, but $Z(M)$
may act on $\sigma$ by an arbitrary character. We call $(M,\sigma )$ a cuspidal 
pair, and from it we construct the parabolically induced $G$-representation 
$I_P^G (\sigma )$. For every $(\pi ,V) \in \mr{Irr} (G)$ there is a cuspidal pair 
$(M, \sigma )$, uniquely determined up to $G$-conjugacy, such that 
$V \in \mr{JH} (I_P^G (\sigma ) )$.

We denote the complex torus of nonramified characters of $M$ by
$\xnr (M)$, and the compact subtorus of unitary nonramified
characters by $\xunr (M)$. We say that two cuspidal pairs
$(M,\sigma )$ and $M' ,\sigma' )$ are inertially equivalent if
there exist $\chi \in \xnr (M')$ and $g \in G$ such that $M' = g
M g^{-1}$ and $\sigma' \otimes \chi \cong \sigma^g$. With an
inertial equivalence class $\mf s = [M,\sigma ]_G$ we associate a
subcategory Rep$ (G)^{\mf s}$ of Rep$ (G)$. By definition its
objects are smooth $G$-representations $\pi$ with the following
property: for every $\rho \in \mr{JH} (\pi)$ there
is a $(M,\sigma ) \in \mf s$ such that $\rho$ is a
subrepresentation of $I_P^G (\sigma )$. These blocks Rep$ (G)^{\mf
s}$ give rise to the Bernstein decomposition \cite[Proposition 2.10]{BeDe}
\[
\mr{Rep} (G) = {\ts \prod_{\mf s \in \Omega (G)} } \, \mr{Rep} (G)^{\mf s} \,.
\]
The set $\Omega (G)$ of Bernstein components is countably infinite. 
There are corresponding decompositions of the Hecke and Schwartz 
algebras of $G$ into two-sided ideals:
\begin{align}
\label{eq:2.3} & \mc H (G) = {\ts \bigoplus_{\mf s \in \Omega (G)} }
\, \mc H (G)^{\mf s} \,, \\
\label{eq:2.4} & \mc S (G) \, = {\ts \bigoplus_{\mf s \in \Omega (G)} }
\, \mc S (G)^{\mf s} \,.
\end{align}
For $C_r^* (G)$ this is a less straightforward, since its elements can be
supported on infinitely many Bernstein components. Let $C_r^* (G)^{\mf s}$ 
be the two-sided ideal generated by $\mc H (G)^{\mf s}$. The reduced 
$C^*$-algebra of $G$ decomposes as a direct sum in the $C^*$-algebra sense:
\begin{equation}\label{eq:2.13}
C_r^* (G) = \varinjlim_{\mf S} \bigoplus_{\mf s \in \mf S} C_r^* (G)^{\mf s} \,,
\end{equation}
where the direct limit runs over all finite subsets $\mf S$ of $\Omega (G)$.
For $K \in \mr{CO}(G)$ and $\mf s \in \Omega (G)$ we write
\[
\begin{array}{ccccc}
\mc H (G,K)^{\mf s} & = & \mc H (G)^{\mf s} & \cap & \mc H (G,K) \,, \\
\mc S (G,K)^{\mf s} & = & \mc S (G)^{\mf s} & \cap & \mc S (G,K) \,, \\
C_r^* (G,K)^{\mf s} & = & C_r^* (G)^{\mf s} & \cap & C_r^* (G,K) \,.
\end{array}
\]
Every element of $\mc H (G)$ has a unique decomposition as a sum
of a part in $\mc H (G)^{\mf s}$ and a part in the annihilator of
this ideal. In particular we can write
\begin{equation}
e_K = e_K^{\mf s} + e'_K \in \mc H (G)^{\mf s} \oplus
\bigoplus_{\mf s' \in \Omega (G) \setminus \{\mf s\} } \mc H (G)^{\mf s'} .
\end{equation}
\begin{prop}\label{prop:2.1}
\begin{description}
\item[a)] For fixed $K \in \mr{CO}(G)$ there exist only finitely
many $\mf s \in \Omega (G)$ such that $\mc H (G,K)^{\mf s} \neq 0$.
\item[b)] For every $\mf s \in \Omega (G)$ there exists a $K_{\mf s}
\in \mr{CO}(G)$ such that for all compact open subgroups $K
\subset K_{\mf s}$ the bimodules $e_K^{\mf s} \mc H (G)$ and $\mc
H (G) e_K^{\mf s}$ provide a Morita equivalence between 
\[
\mc H (G)^{\mf s} = \mc H (G) e_K^{\mf s} \mc H (G) \quad \text{and}  
\quad \mc H (G,K)^{\mf s} = e_K^{\mf s} \mc H (G) e_K^{\mf s} \,.
\]
\item[c)] As \textup{b)}, but with $\mc S (G)$ instead of $\mc H (G)$.
\item[d)] As \textup{b)}, but with $C_r^* (G)$ instead of $\mc H (G)$.
\end{description}
\end{prop}
\emph{Proof.} a) See \cite[\S 3.7]{BeDe}.\\
b) By \cite[Corollaire 3.9]{BeDe} there exists a $K_{\mf s} \in \mathrm{CO}(G)$ such that
\[
\mc H (G)^{\mf s} = \mc H (G) e_{K_{\mf s}}^{\mf s} \mc H (G) \,.
\]
For any $K \in \mr{CO}(G)$ with $K \subset K_{\mf s}$ we have 
$e_{K_{\mf s}} \in \mc H (G) e_K \mc H (G)$. Hence
\[
\mc H (G) e_K^{\mf s} \mc H (G) = \big( \mc H (G) e_K \mc H (G) \big) \cap \mc H (G)^{\mf s}
\]
contains $e_{K_{\mf s}}^{\mf s}$ and must equal $\mc H (G)^{\mf s}$.
We note that these $e^{\mf s}_K$ are precisely the special idempotents
constructed in \cite[Proposition 3.13]{BuKu}. \\
c) and d) follow directly from b) and the characterisation of 
$\mc S (G)^{\mf s}$ and $C_r^* (G)^{\mf s}$ as the ideals 
generated by $\mc H (G)^{\mf s}. \qquad \Box$
\\[2mm]

Let $\mf s = [M,\sigma ]_G \in \Omega (G)$ and put
\[
\begin{array}{lll}
N (M,\sigma ) & = & \{ g \in N_G (M) : \sigma^g \cong \sigma \otimes \chi 
\text{ for some } \chi \in \xnr (M) \} \,, \\
W_{\mf s} & = & N (M,\sigma ) / M \,.
\end{array}
\]
\begin{thm}\label{thm:2.15}
For $K \in \mr{CO}(G)$ as in Proposition \textup{\ref{prop:2.1}.b} $\mc H (G,K)^{\mf s}$ 
is a unital finite type algebra with center isomorphic to 
\[
\mc O (\xnr (M) / W_{\mf s} ) = \mc O (\xnr (M) )^{W_{\mf s}} .
\]
There are natural isomorphisms 
\[
\begin{array}{lll}
Z \big( \mc H(G,K_{\mf s})^{\mf s} \big) & \cong & Z \big( \mr{Rep}(G )^{\mf s} \big) \,, \\
\mr{Prim} \big( \mc H(G,K_{\mf s})^{\mf s} \big) & \cong & 
 \mr{Irr}(G) \cap \mr{Rep}(G )^{\mf s} \,, \\
\mr{Prim} \big( \mc S(G,K_{\mf s})^{\mf s} \big) & \cong &
 \mr{Irr}^t(G) \cap \mr{Rep}(G )^{\mf s} \,.
\end{array}
\]
\end{thm}
\emph{Proof.}
This follows from Th\'eor\`eme 2.13 and Corollaire 3.4 of \cite{BeDe},
in combination with Proposition \ref{prop:2.1}. $\qquad \Box$
\\[2mm]

The cohomological dimension of the abelian category Rep$ (G)$ equals the 
rank of $G$, which is by definition the dimension of a maximal split
subtorus of $G$ \cite[Section II.3]{ScSt}. The author does not
know whether the abelian category $\mr{Rep}^t (G)$ has finite
cohomological dimension. Yet one can determine something like
the global dimension of $\mc S (G)$, at the cost of using more
advanced techniques. Namely, according to Meyer \cite[Theorem 29]{Mey} 
the cohomological dimension of the exact category
$\mr{Mod}_b (\mc S (G))$ of complete essential bornological
$\mc S (G)$-modules is also equal to the rank of $G$. The natural
tensor product to work with in this category is the completed
bornological $A$-balanced tensor product, which we denote by
$\hot_A$. For Fr\'echet spaces $\hot_{\mh C}$ agrees with the
completed projective tensor product, so we abbreviate it to $\hot$.
For later use we translate these cohomological dimensions to 
statements about Hochschild homology.

\begin{lem}\label{lem:2.2}
\begin{description}
\item[a)] $HH_n (\mc H (G)) = 0$ for all $n > \mr{rk}(G) $,
\item[b)] $HH_n (\mc S (G),\hot_{\mh C} ) = 0$ for all $n > \mr{rk}(G) $.
\end{description}
\end{lem}
\emph{Proof.} a) can be found in \cite{Nis} but we prefer to
derive it from the above. Let
\begin{equation}\label{eq:2.15}
\mh C \leftarrow P_0 \leftarrow P_1 \leftarrow \cdots \leftarrow
P_{\mr{rk}(G)} \leftarrow 0
\end{equation}
be a projective resolution of the trivial $G$-module $\mh C$.
Endowing $P_m \otimes \mc H (G)$ with the diagonal $G$-action,
$\mc H (G) \leftarrow P_* \otimes \mc H (G)$ becomes a resolution
of $\mc H (G)$ by projective $\mc H (G)$-bimodules. By definition
\begin{equation}
\begin{split}
HH_n (\mc H (G)) & = \mr{Tor}_n^{\mc H (G) \otimes \mc H (G)^\mr{op}}
(\mc H (G), \mc H (G)) \\
& = H_n \big( \mr{Hom}_{\mc H (G) \otimes \mc H (G)^\mr{op}} \big(
P_* \otimes \mc H (G) , \mc H (G) \big) \big) \,,
\end{split}
\end{equation}
which clearly vanishes for $n > \mr{rk}(G)$.\\
b) We will use that the inclusion $\mc H (G) \to \mc S (G)$ is
isocohomological \cite[Theorem 22]{Mey}. According to
\cite[(22)]{Mey} the differential complex
\[
\mc S (G) \leftarrow \mc S (G) \hot_{\mc H (G)} P_* \hot_{\mh C} \mc S (G)
\]
is a projective resolution of $\mc S (G)$ in $\mr{Mod}_b (\mc S (G))$. Hence
\begin{equation}\label{eq:2.5}
\begin{split}
HH_n (\mc S (G), \hot_{\mh C} ) & = \mr{Tor}_n^{\mc S (G)
\hot \mc S (G)^\mr{op}} (\mc S (G),\mc S (G)) \\
& = H_n \big( \mr{Hom}_{\mc S (G) \hot \mc S (G)^\mr{op}} \big( \mc S (G) 
\hot_{\mc H (G)} P_* \hot_{\mh C} \mc S (G) , \mc S (G) \big) \big) \,.
\end{split}
\end{equation}
From \eqref{eq:2.15} we see immediately that this vanishes 
$\forall n > \mr{rk}(G). \qquad \Box$
\\[3mm]
\begin{cor}\label{cor:2.3}
Let $\mf s \in \Omega (G)\,, n > \mr{rk}(G)$ and $K \in \mr{CO}(G)$
be such that $K \subset K_{\mf s}$.
\begin{description}
\item[a)] $HH_n (\mc H (G)^{\mf s}) = 0 = HH_n (\mc H (G,K)^{\mf s})$ ,
\item[b)] $HH_n (\mc S (G)^{\mf s},\hot_{\mh C} ) = 0 =
HH_n (\mc S (G,K)^{\mf s},\hot)$ .
\end{description}
\end{cor}
\emph{Proof.}
b) From \eqref{eq:2.4} and \eqref{eq:2.5} we see that
\begin{equation}
HH_n (\mc S (G),\hot_{\mh C} ) \cong {\ts \bigoplus_{\mf s \in \Omega (G)} }
HH_n (\mc S (G)^{\mf s},\hot_{\mh C} ) \,.
\end{equation}
By Proposition \ref{prop:2.1}.c
\begin{equation}
\begin{split}
& HH_n (\mc S (G,K)^{\mf s},\hot) = \mr{Tor}_n^{\mc S (G,K)^{\mf s} \hot 
\mc S (G,K)^{\mf s,\mr{op}}} (\mc S (G,K)^{\mf s}, \mc S (G,K)^{\mf s}) \cong \\
& \mr{Tor}_n^{\mc S (G)^{\mf s} \hot \mc S (G)^{\mf s,\mr{op}}}
(\mc S (G)^{\mf s}, \mc S (G)^{\mf s}) = HH_n (\mc S (G)^{\mf s},\hot_{\mh C} ) \,,
\end{split}
\end{equation}
where we take the torsion functors in the category of complete
bornological modules. By Lemma \ref{lem:2.2}.b these homology
groups all vanish for $n > \mr{rk}(G)$.\\
a) can be proved in exactly the same way as b), using \eqref{eq:2.3}, 
Lemma \ref{lem:2.2}.a and Proposition \ref{prop:2.1}.b. $ \qquad \Box$ 
\vspace{4mm}

\section{The Plancherel theorem}
\label{sec:2.2}

The Plancherel formula for $G$ is an explicit decomposition of the trace
\[
\mc H (G) \to \mh C \quad,\quad f \mapsto f (e)
\]
in terms of the traces of irreducible $G$-representations. Closely
related is the Plan\-che\-rel theorem, which describes $\mc S (G)$ in
terms of its irreducible representations. This description is due
to Harish-Chandra \cite{HC1,HC2}, although he published only a
sketch of the proof. Harish-Chandra's notes were worked out in
detail by Waldspurger \cite{Wal}. In the present section we recall the
most important ingredients of the Plancherel theorem, relying
almost entirely on the above papers.

A parabolic pair $(P,A)$ consists of a parabolic subgroup $P$ of
$G$ and a maximal split torus $A$ in
the center of some Levi subgroup $M$ of $P$. If $N$ is the
unipotent radical of $P$ then $P = M \ltimes N$ and $M = Z_G (A)$.
Moreover restriction from $M$ to $A$ defines a surjection
$\xnr (M) \to \xnr (A)$ with finite kernel.

The maximal parabolic pair is $(G,A_G )$, where $A_G$ is the unique
maximal split torus of $Z(G)$.
We fix a maximal split torus $A_0$ of $G$, and a minimal parabolic
subgroup $P_0$ containing $A_0$. We call $(P,A)$ semi-standard if
$A \subset A_0$, and standard if moreover $P \supset P_0$. Every
parabolic pair is conjugate to a standard one.

Let $(\omega ,E)$ be an irreducible square-integrable representation of $M$. 
By definition this entails that $E$ is smooth, pre-unitary and admissible. 
Let $(\breve \omega, \breve E)$ be the smooth contragredient representation. 
The admissible $G\times G$-representation
\[
L(\omega ,P) = I_{P \times P}^{G \times G} (E \otimes \breve E) =
I_P^G (E) \otimes I_P^G (\breve E )
\]
is naturally a nonunital Hilbert algebra. Notice that for every
$\chi \in \xunr (M)$ the representation $\omega \otimes \chi$ is
also square-integrable, and that $L(\omega \otimes \chi ,P)$ can
be identified with $L(\omega ,P)$. Let $\mb k_\omega$ be the set of $k
\in \xnr (M)$ such that $\omega \otimes k$ is equivalent to
$\omega$. This is a finite subgroup of $\xunr (M)$. For every $k
\in \mb k_\omega$ there exists a canonical unitary intertwiner
\[
I(k,\omega ) \in \mr{Hom}_{G \times G} (L (\omega ,P),L(\omega \otimes k,P)) \,.
\]

Next we consider the intertwiners associated to elements of
various Weyl groups. Let  $(Q,B)$ be another parabolic pair. Write
$W(A|G|B)$ for the set of all homomorphisms $B \to A$ induced by
inner automorphisms of $G$. This is a group in case $A = B $:
\[
W(G,A) := W(A|G|A) = N_G (A) / Z_G (A) = N_G (A) / M \,.
\]
Let $(Q,A^g )$, with $g \in G$, be yet another parabolic pair, and
put $n = [g] \in W (A^g |G|A)$. The equivalence class of the
$M^g$-representation $(\omega^{g^{-1}},E)$ depends only on $n$ and
may therefore be denoted by $n \omega$. Waldspurger constructs
certain normalized intertwiners $\prefix{^\circ}{c_{Q|P}}(n,\omega
)$. Preferring the simpler notation $I(n,\omega)$ we recall their
properties.

\begin{thm}\label{thm:2.4} \textup{\cite[Paragraphe V]{Wal}} \\
Let $(P,A) ,\, (P',A')$ and $(Q,B)$ be semi-standard
p-pairs, and $n \in W (B|G|A)$. There exists an intertwiner
\[
I (n,\omega \otimes \chi) \in \mr{Hom}_{G \times G} (L (\omega ,P)
, L (n \omega ,Q) )
\]
with the following properties:
\begin{itemize}
\item $\chi \to I (n, \omega \otimes \chi)$ is a rational
function on $\xnr (M) $,
\item $I (n ,\omega \otimes \chi)$ is unitary and regular
for $\chi \in \xunr (M) $,
\item If $n' \in W (A' |G|B)$ then
\[
I (n',n (\omega \otimes \chi) ) I (n,\omega \otimes \chi) = 
I (n' \, n, \omega \otimes \chi) \,.
\]
\end{itemize}
\end{thm}

To define the Fourier transform implementing the Plancherel
isomorphism we introduce a space of induction data, such that
every irreducible tempered representation is a direct summand of
(at least) one of these parabolically induced representations. For
every semi-standard parabolic pair $(P,A)$ choose a set 
$\Delta_M$ of irreducible square-integrable representations of
$M = Z_G (A)$, with the following property. For every
square-integrable $\pi \in \mr{Irr}(M)$ there exists precisely one
$\omega \in \Delta_M$ such that $\pi$ is equivalent to $\omega \otimes \chi$, 
for some $\chi \in \xnr (M)$. We call a triple $(P,A,\omega )$ standard if 
$(P,A)$ is a standard parabolic pair and $\omega \in \Delta_M$.

An induction datum is a quadruple $(P,A,\omega,\chi )$ where
$(P,A)$ is a semi-standard parabolic pair, $\omega \in \Delta_M$
and $\chi \in \xnr (M)$. Let $\Xi$ be the scheme of all induction
data and $\Xi_u$ the smooth submanifold of unitary induction data,
that is, those with $\chi \in \xunr (M)$. Then $\Xi$ and $\Xi_u$ are
countable disjoint unions of complex algebraic tori and compact tori,
respectively. For $\xi = (P,A,\omega,\chi ) \in \Xi$ we put
\[
I(\xi) = I_P^G (\omega \otimes \chi) \,.
\]
By \cite[Lemme III.2.3]{Wal} the representation $I(\xi)$ is
tempered if and only if $\omega \otimes \chi$ is tempered, if and
only if $\xi \in \Xi_u$. Like for cuspidal pairs one can define
inertial equivalence on $\Xi_u$. The set $\Omega^t (G)$ of all
equivalence classes $[P,A,\omega]_G$ is called the Harish-Chandra
spectrum of $G$. It comes with a natural surjection
$\Omega^t (G) \to \Omega (G)$, see \cite[Section 1]{SSZ}. It follows
from Proposition \ref{prop:2.1} and \cite[Th\'eor\`eme VIII.1.2]{Wal} 
that this map is finite-to-one.

Let $\mc L_{\Xi}$ be the vector bundle over $\Xi$ which is trivial
on every component and whose fiber at $\xi$ is $L(\omega,P)$. We
say that a section of this bundle is algebraic (polynomial) or
rational if it is supported on only finitely many components, and
has the required property on every component. Now we can define
the Fourier transform:
\begin{equation}\label{eq:2.6}
\begin{aligned}
& \mc F : \mc H (G) \to \mc O (\Xi ; \mc L_{\Xi}) \,, \\
& \mc F (f) (P,A,\omega,\chi ) = I(P,A,\omega,\chi)(f) \in
L(\omega ,P) \,,
\end{aligned}
\end{equation}
where we used the notation from \eqref{eq:2.12}.
Notice that this differs slightly from $\check f (\omega \otimes
\chi,P)$ as in \cite[\S VII.1]{Wal}. To make it fit better with
its natural adjoint Waldspurger adjusts the Fourier transform. We
will use \eqref{eq:2.6} though, because it is multiplicative.

To formalize the action of the intertwiners on sections of $\mc
L_\Xi$ we construct a locally finite groupoid $\mc W$. The set of
objects of $\mc W$ is $\Xi$ and the morphisms from $\xi$ to $\xi'$
are the pairs $(k,n)$ with the following properties:
\begin{itemize}
\item $k \in \mb k_\omega$ ,
\item $n \in W(A|G|A')$ and $n (A') = A$,
\item $n \omega'$ is equivalent to $\omega \otimes \tilde \chi$
for some $\tilde \chi \in \xnr (M) $.
\end{itemize}
The multiplication in $\mc W$, if possible, is
\[
(k,n) (k',n') = (k (k' \circ n),n n') \,.
\]
Let $\Gamma (\Xi ;\mc L_\Xi)$ be a suitable algebra of sections of
$\mc L_\Xi$. For $f \in \Gamma (\Xi ;\mc L_\Xi)$ we define
\begin{align*}
& k \cdot f (\chi) = I(k, \omega) f(k^{-1} \omega) \,,\\
& n \cdot f (\chi' ) = I(n,\omega ) f (\chi' \circ n) \,.
\end{align*}
Notice that the intertwiners do not stabilize $\mc O (\Xi ;\mc L_\Xi )$ in general.
Nevertheless, we write $\mc O (\Xi ;\mc L_\Xi )^{\mc W}$ for the subalgebra of
$\mc W$-invariant sections, which by construction contains $\mc F (\mc H (G))$.
Because $(\omega,E)$ is admissible,
\[
C^\infty (\xunr (M)) \otimes L(\omega,P)^{K \times K} = C^\infty
(\xunr (M)) \otimes I_P^G (E)^K \otimes I_P^G (\breve E )^K
\]
has a natural Fr\'echet topology, for every $K \in \mr{CO}(G)$. We endow
\[
C_c^\infty (\Xi_u ;\mc L_\Xi) = \varinjlim_{K \in
\mr{CO}(G)} C_c^\infty \big( \Xi_u ; \mc L_\Xi^{K \times K} \big)
\]
with the inductive limit topology. The Plancherel theorem for
reductive $p$-adic groups reads:

\begin{thm}\label{thm:2.5} \textup{\cite{HC2,Wal}} \\
The Fourier transform
\[
\mc F : \mc S (G) \to C_c^\infty (\Xi_u ;\mc L_\Xi )^{\mc W}
\]
is an isomorphism of topological algebras.
\end{thm}

A simple representation theoretic consequence of this important theorem is

\begin{cor}\label{cor:2.15}
For any $w \in \mc W$ and $\xi \in \Xi$ such that $w \xi$ is
defined, the $G$-representations $I(\xi )$ and $I(w \xi )$ have the
same irreducible subquotients, counted with multiplicity.
\end{cor}
\emph{Proof.}
By \cite[Corollary 2.3.3]{Cas} we have to show that the traces of $I(\xi )$ and 
$I(w \xi )$ are the same, in other words, that the function
\begin{equation*}
\mc H (G) \times \xnr (M) \to \mh C : (f ,\chi) \mapsto
\mr{tr}\, I(P,A,\omega ,\chi )(f) - 
\mr{tr}\, I (wP,wA,w \omega ,\chi \circ w^{-1})(f)
\end{equation*}
is identically 0. Because this is a polynomial function of $\chi$,
it suffices to show that it is 0 on $\mc H (G) \times \xunr (M)$.
That follows from Theorem \ref{thm:2.5}. $\qquad \Box$
\\[2mm]

The Plancherel theorem can be used to describe the Fourier transform of
$C_r^* (G)$. For $(\omega,E) \in \Delta_M$ let $\mc K (\omega ,P)$
be the algebra of compact operators on the Hilbert space
completion of $I_P^G (E)$. Notice that
\begin{equation}
\mc K (\omega ,P) = \varinjlim_{K \in \mr{CO}(G)} L (\omega ,P)^{K
\times K}
\end{equation}
in the $C^*$-algebra sense, and that the intertwiner $I (n,
\omega)$ extends to $\mc K (\omega ,P)$ because it is unitary. Let
$\mc K_\Xi$ be the vector bundle over $\Xi$ whose fiber at
$(P,A,\omega ,\chi )$ is $\mc K (\omega ,P)$, and let $C_0 (\Xi_u ;
\mc K_\Xi )$ be the $C^*$-completion of
\[
{\ts \bigoplus_{(P,A,\omega)} } C (\xunr (M) ; \mc K (\omega ,P)) \,.
\]
\begin{thm}\label{thm:2.6} \textup{\cite[Theorem 2.5]{Ply}} \\
The Fourier transform extends to an isomorphism of $C^*$-algebras
\[
C_r^* (G) \to C_0 (\Xi_u ; \mc K_\Xi )^\mc W \,.
\]
\end{thm}

We can also describe the images of the subalgebras $\mc S (G,K)$
and $C_r^* (G,K)$ under the Fourier transform.

\begin{thm}\label{thm:2.7}
Fix $K \in \mr{CO}(G)$. There exists a finite set of standard triples 
$(P_i ,A_i ,\omega_i )$ with the following properties.
\begin{description}
\item[a)] The Fourier transform yields algebra homomorphisms
\[
\begin{array}{rrr}
\mc H (G,K) & \to & \bigoplus_{i=1}^{n_K}
  \big( \mc O (\xnr (M_i )) \otimes L(\omega_i ,P_i
  )^{K \times K} \big)^{\mc W_i} \\
\mc S (G,K) & \to & \bigoplus_{i=1}^{n_K}
  \big( C^\infty (\xunr (M_i )) \otimes L(\omega_i ,P_i
  )^{K \times K} \big)^{\mc W_i} \\
C_r^* (G,K) & \to & \bigoplus_{i=1}^{n_K}
  \big( C (\xunr (M_i )) \otimes L(\omega_i ,P_i )^{K \times K} \big)^{\mc W_i}
\end{array}
\]
where $\mc W_i$ is the isotropy group of $(P_i ,A_i ,\omega_i )$
in $\mc W$.
\item[b)] The first map is injective, the second is an isomorphism of Fr\'echet 
algebras and the third is an isomorphism of $C^*$-algebras.
\item[c)] For every $w \in \mc W_i$ there exists a unitary section
\[
u_w \in C^\infty (\xunr (M_i )) \otimes L(\omega_i ,P_i )^{K
\times K}
\]
which extends to a rational section on $\xnr (M_i )$, such that
\end{description}
\begin{equation}\label{eq:2.7}
w f (\chi ) = u_w (\chi ) f (w^{-1} \chi ) u_w^{-1}(\chi ) \qquad
\forall f \in C (\xunr (M_i )) \otimes L(\omega_i ,P_i )^{K
\times K} .
\end{equation}
\end{thm}

\emph{Proof.} The author already proved this result in
\cite[Theorem 10]{Sol1} but we include the proof anyway. Notice
that, in constrast with Proposition \ref{prop:2.1}, it is not necessary
to require that $K$ is ``small", expect for being compact. That is
because the most tricky (namely, not completely reducible)
representations in a Bernstein component $[M,\sigma ]_G$ appear
only if we twist $\sigma$ by a nonunitary character $\chi \in
\xnr (M) \setminus \xunr (M)$.

According to \cite[Th\'eor\`eme VIII.1.2]{Wal} there are only
finitely many components in the Harish-Chandra spectrum $\Omega^t (G)$
on which the idempotent $e_K$ does not act as 0. Pick one triple
$(P_i ,A_i ,\omega_i )$ for each such component. Now a) and b)
follow immediately from Theorems \ref{thm:2.5} and \ref{thm:2.6}.

Concerning c), every automorphism of
\[
L(\omega_i ,P_i )^{K \times K} \cong \mr{End}_{\mh C} \big( I_P^G
(E)^K \big)
\]
is inner, so \eqref{eq:2.7} holds for some section $u_w$. Using
Theorem \ref{thm:2.4} we can arrange that $u_w$ is rational on
$\xnr (M_i )$ and unitary on $\xunr (M_i ). \qquad \Box$
\vspace{4mm}

\section{The Langlands classification}
\label{sec:2.3}

The Langlands classification describes the relation between the
smooth spectrum of $G$ and the tempered spectra of its Levi subgroups.
Let $(P,A)$ be a semi-standard parabolic pair, let $X^* (A)$ be the lattice of 
algebraic characters of $A$ and put 
\[
\mf a^* = X^* (A) \otimes_{\mh Z} \mh R \cong X^* (M) \otimes_{\mh Z} \mh R \,.
\]
For $A = A_0$ and $A = A_G$ we write $\mf a^* = \mf a_0^*$ and 
$\mf a^* = \mf a_G^*$, respectively. There is a natural homomorphism
\begin{equation}
\begin{split}
& H_M : M \to \mathrm{Hom}_{\mh Z} (X^* (M), \mh R) \cong
  \mathrm{Hom}_{\mh R} (\mf a^* ,\mh R ) \,, \\
& q^{\inp{\chi}{H_M (m)}} = \norm{\chi (m)}_{\mh F} 
\qquad \qquad m \in M , \chi \in X^* (M) \,,
\end{split}
\end{equation}
where $q$ is the cardinality of the residue field of $\mh F$. Conversely, for 
$\nu \in \mf a^*$ we define a nonramified character $\chi_\nu$ of $M$ by
\begin{equation}
\chi_\nu (m) = q^{\inp{\nu}{H_M (m)}} .
\end{equation}
This yields an isomorphism 
\[
\mf a^* \cong \mr{Hom}(M, \mh R_{>0}) \,,
\]
where Hom is taken in the category of topological groups.
Let $Q$ be a parabolic subgroup such that $P \subset Q \subset G$,
and let $\mf q$ be its Lie algebra. It decomposes into $A$-eigenspaces
\begin{equation}
\mf q_\alpha := \{ x \in \mf q : \mr{Ad}(a) x = \alpha (a) x \;
\forall a \in A \}
\end{equation}
with $\alpha \in X^* (A)$. The roots of $Q$ with respect to $A$ are
\begin{equation}
\Sigma (Q,A) := \{ \alpha \in X^* (A) \setminus \{ 1 \} : \mf
q_\alpha \neq 0 \} \,.
\end{equation}
The inclusions $A \to A_0$ and $M_0 \to M$ identify $\mf a^*$ as a direct 
summand of $\mf a_0^*$. We have a root system
\[
\Sigma_0 = \Sigma (G,A_0 ) \subset \mf a_0^*
\]
with positive roots $\Sigma (P_0 ,A_0 )$, simple roots $\Delta_0 = \Delta(P_0, A_0)$ 
and Weyl group $W_0 = W (P_0 ,A_0 )$. 
We fix a $W_0$-invariant inner product $\langle \,, \rangle_0$ 
on $\mf a_0^*$, so that we may identify this vector space with its dual. 
Denote the Lie algebras of $A$ and $A_0$ by $\mf a$ and $\mf a_0$. 
(Notice that $\mf a_0$ is not the dual of $\mf a_0^* $, these are vector 
spaces over different fields.) For $F \subset \Delta_0$ we put
\begin{equation}\label{eq:2.8}
\begin{array}{lll@{\qquad}lll}
\Sigma_F & = & \Sigma_0 \cap \mh R F \; \subset \; \mf a_0^* &
W_F & = & \langle s_\alpha : \alpha \in F \rangle \; \subset \; W_0 \,, \\
\mf a_F & = & \{ x \in \mf a_0 : \alpha (x) = 0 \, \forall
\alpha \in F \} & A_F & = & \exp (\mf a_F ) \,, \\
\mf m_F & = & \big( \bigoplus_{\alpha \in \Sigma_F} \mf g_\alpha
\big) \oplus \mf a_0 & M_F & = & Z_G (A_F ) \,, \\
\mf n_F & = & \bigoplus_{\alpha \in \Sigma (P_0 ,A_0 ) \setminus
\Sigma_F} \mf g_\alpha & N_F & = & \exp (\mf n_F ) \,, \\
\mf p_F & = & \big( \bigoplus_{\alpha \in \Sigma (P_0 ,A_0 ) \cup
\Sigma_F} \mf g_\alpha \big) \oplus \mf a_0 & P_F & = & M_F
\ltimes N_F \,.
\end{array}
\end{equation}
Every standard parabolic pair is of the form $(P_F ,A_F )$ for
some $F \subset \Delta_0$. In this situation $F = \Delta (P_F ,A_F )$ 
is the set of nonzero projections of $\Delta_0 \subset \mf a_0^*$ on 
$\mf a_F^*$, and $W(M_F ,A_0 ) = W_F$. In particular a standard parabolic 
pair is completely determined by either of its two ingredients.

Furthermore we introduce the open and closed positive cones in $\mf a^*$:
\begin{align*}
& \mf a^{*,+} = \{ \nu \in \mf a^* : \inp{\nu}{\alpha}_0 > 0
\,\forall \alpha \in \Delta (P,A) \} \,, \\
& \bar{\mf a}^{*,+} = \{ \nu \in \mf a^* : \inp{\nu}{\alpha}_0 \geq
0 \,\forall \alpha \in \Delta (P,A) \} \,.
\end{align*}
Their antidual is the obtuse negative cone in $\mf a_0^*$ : 
%& \mf a_0^{*,-} = \{ \mu \in \mf a_0^* : \inp{\nu}{\mu} < 0 \; 
%\forall \nu \in \mf a_0^{*,+} \} \,, \\
\[
\bar{\mf a}_0^{*,-} := \{ \mu \in \mf a_0^* : \inp{\nu}{\mu} \leq 0 \; 
\forall \nu \in  \bar{\mf a}_0^{*,+} \} \,.
\]
The set of Langlands data $\Lambda^+$ consists of all quadruples
$\lambda = (P,A,\sigma,\nu )$ such that
\begin{itemize}
\item $(P,A)$ is a standard parabolic pair with Levi component $M
= Z_G (A) $,
\item $\sigma \in \mr{Irr}^t (M) $,
\item $\nu \in \mf a^{*,+} $.
\end{itemize}
Given a Langlands datum $\lambda \in \Lambda^+$ we pick a concrete realization of
$\sigma$ and we construct the admissible $G$-representation 
$I (\lambda ) = I_P^G (\sigma \otimes \chi_\nu )$. Notice that $\lambda$ determines
$I(\lambda )$ only modulo equivalence of $G$-representations. 
The classical Langlands classification for reductive $p$-adic groups reads:

\begin{thm}\label{thm:2.8}
\textup{\textbf{a)}} For every $\lambda \in \Lambda^+$ the $G$-representation
$I (\lambda)$ is indecomposable and has a unique irreducible quotient.\\
We call this the Langlands quotient $J (\lambda )$.\\
\textup{\textbf{b)}} For every $\pi \in \mr{Irr}(G)$ there is a unique
$\lambda \in \Lambda^+$ such that $\pi$ is equivalent to $J (\lambda )$.
\end{thm}
\emph{Proof.} See \cite{Kon} or \cite[\S XI.2]{BoWa}. We note that Konno proves
the uniqueness part only modulo $W_0$-conjugacy. However, two $W_0$-conjugate
Langlands data are necessarily equal, which we will show in a more general setting
in Lemma \ref{lem:2.10}. $\qquad \Box$ 
\\[3mm]

The Langlands datum of $\pi \in \mr{Irr}^t (G)$ is simply $(G,A_G ,\pi ,0)$, and conversely
$J(\lambda )$ cannot be tempered if $\nu \neq 0$. For every $\lambda \in \Lambda^+$
the Langlands quotient is maximal among the irreducible constituents of $I(\lambda )$,
in a suitable sense:

\begin{lem}\label{lem:2.9}
Let $\lambda = (P,A,\sigma,\nu)$ and $\lambda' = (P',A',\sigma',\nu')$ be Langlands data,
and suppose that $J(\lambda' ) \in \mr{JH}(I(\lambda ))$.
\begin{description}
\item[a)] $\nu' - \nu \in \bar{\mf a}_0^{*,-} , A' \subset A$ and $P' \supset P$.
\item[b)] If $\nu' = \nu$ then $\lambda' = \lambda$.
\item[c)] $\mr{End}_G (I(\lambda )) = \mh C$.
\end{description}
\end{lem}
\emph{Proof.} a) and b) The statements about $\nu$ and $\nu'$ are 
\cite[Lemma XI.2.13]{BoWa}. From the definition of $\Lambda^+$ we see that $A' \subset A$ 
and $P' \supset P$ whenever $\nu' - \nu \in \bar{\mf a}_0^{*,-}$.\\
c) Let $\phi \in \mr{End}_G (I(\lambda ))$ and write $M(\lambda ) = \ker (I (\lambda) \to 
J (\lambda ))$. By Theorem \ref{thm:2.8} and the above $M(\lambda)$ and $J(\lambda )$ 
do not have any common irreducible constituents. Hence the composition 
\[
M (\lambda) \xrightarrow{\phi} I(\lambda ) \to J(\lambda )
\]
of $\phi$ with the quotient map is zero, and $\phi (M(\lambda )) \subset M(\lambda)$. 
Therefore $\phi$ induces a map $\phi_J \in \mr{End}_G (J (\lambda))$ and, because 
$J(\lambda )$ is irreducible, $\phi_J = \mu \, \mr{Id}_{J(\lambda )}$ for some 
$\mu \in \mh C$. We want to show that
\[
\psi := \phi - \mu \, \mr{Id}_{I(\lambda )} \in \mr{End}_G (I(\lambda ))
\]
equals zero. By construction $\psi (I(\lambda )) \subset M (\lambda )$. Let $V$ be an 
irreducible quotient representation of $M(\lambda )$ and consider the induced map 
$\psi_V \in \mr{Hom}_G (I(\lambda ),V)$. It follows from Theorem \ref{thm:2.8}.a that 
the smooth contragredient representation $\check{I (\lambda )}$ of $I(\lambda )$ has 
exactly one irreducible submodule, which moreover is equivalent to $\check{J(\lambda )}$. 
Since $\check V$ is not equivalent to the contragredient of $J(\lambda )$, we find
\[
0 = \mr{Hom}_G (\check V , \check{I (\lambda )}) \cong \mr{Hom}_G (I (\lambda ),V) \,.
\]
In particular $\psi_V = 0$ and $\psi (I(\lambda)) \subset \ker (M(\lambda) \to V)$. 
Since $M(\lambda )$ has finite length \cite[p. 30]{BeDe}, we conclude with induction that 
\begin{equation}\label{eq:2.19}
\mr{Hom}_G (I(\lambda ),M(\lambda )) = 0 \,.
\end{equation}
Thus $\psi = 0$ and $\phi = \mu \, \mr{Id}_{I(\lambda )} \,. \qquad \Box$
\\[2mm]

We would like to reformulate Lemma \ref{lem:2.9}.b with a condition on $\sigma$ instead
of on $\nu$. To achieve this we will define a variation on the central character of a 
representation. Suppose that $\pi \in \mr{Irr}(G)$ belongs to the Bernstein component 
$\mf s = [M ,\rho ]_G$. We may assume that $M$ is a standard Levi subgroup and that $\rho$ 
is a unitary supercuspidal $M$-representation. Pick $\chi_\pi \in \xnr (M)$ such that 
$\pi \in \mr{JH}(I_P^G (\rho \otimes \chi_\pi ))$, and consider $\log |\chi_\pi | \in \mf a^*$. 
This does not depend on the choice of $\rho$, and by Theorem \ref{thm:2.15} it is unique modulo 
$W_{\mf s}$. However, we could have chosen another standard Levi subgroup $M'$ conjugate to 
$M$. Since $W(A' |G| A) \subset W_0$, this would lead to 
\[
\log |\chi'_\pi | = w \log | \chi_\pi | \in \mf a_0^* \quad \text{for some } w \in W_0 \,.
\]
Thus we get an invariant
\[
\cc_G (\pi ) := W_0 \log |\chi_\pi | \in \mf a_0^* / W_0 \,,
\]
which can be considered as a substitute for the absolute value of the $Z(\mc H (G))$-character 
of $\pi$. We note that, because the inner product on $\mf a_0^*$ is $W_0$-invariant, 
$\norm{\cc_G (\pi )}$ is well-defined. Furthermore the orthogonal projection of $\cc_G (\pi )$ 
on $\mf a_G^*$ consists of a single element, known as the central exponent of $\pi$. If 
$\pi \in \mr{JH}(I_{P_F}^G (\tau ))$ for some $\tau \in \mr{Irr}(M_F)$ then, by the transitivity 
of parabolic induction
\begin{equation}\label{eq:2.21}
\cc_G (\pi ) = W_0 \, \cc_{M_F} (\tau ) \,.
\end{equation}

\begin{lem}\label{lem:2.13}
Let $\lambda = (P,A,\sigma,\nu)$ and $\lambda' = (P',A',\sigma',\nu')$ be different 
Langlands data, and suppose that $J(\lambda' ) \in \mr{JH}(I(\lambda ))$. Then 
\[
\norm{\cc_{M'}(\sigma' )} > \norm{\cc_M (\sigma )} \,.
\]
\end{lem}
\emph{Proof.}
By \eqref{eq:2.21} all irreducible constituents of $I(\lambda )$ have the same 
$\cc_G$-invariant, so $\cc_G (J (\lambda' )) = \cc_G (J(\lambda ))$. We have
\[
\cc_M (\sigma \otimes \chi_\nu ) = \cc_M (\sigma ) + \nu = 
W_M \log |\chi_\sigma | + \nu \in \mf a_0^* / W_M \,.
\]
Since $\sigma$ is irreducible and tempered, it is a unitary $M$-representation 
\cite[Proposition III.4.1]{Wal}, in particular $\log |\chi_\sigma | = 0$ on $A$. 
But $\nu \in \mf a^*$ is zero on the derived group of $M$, 
so $\inp{\log |\chi_\sigma |}{\nu}_0 = 0$. Hence
\[
\norm{\cc_G (J(\lambda ))}^2 = \norm{\cc_M (\sigma )}^2 + \norm{\nu}^2 ,
\]
and similarly
\[
\norm{\cc_G (J(\lambda' ))}^2 = \norm{\cc_{M'} (\sigma' )}^2 + \norm{\nu'}^2 .
\]
By Lemma 2.9 $\nu' - \nu \in \bar{\mf a}_0^{*,-} \setminus \{ 0 \}$, which by 
\cite[Claim 3.5.1]{Kon} implies $\norm{\nu} > \norm{\nu'}$. 
Since $\norm{\cc_G (J(\lambda ))}^2 = \norm{\cc_G (J(\lambda' ))}^2$, we conclude that 
$\norm{\cc_M (\sigma )}^2 < \norm{\cc_{M'} (\sigma' )}^2 . \; \Box$
\vspace{4mm}

\section{Parametrizing irreducible representations}
\label{sec:2.4}

We will combine the Langlands classification with the Plancherel theorem to 
parametrize of the irreducible smooth $G$-representations. This parametrization is not 
complete, in the sense that certain packets contain more than one irreducible representation, 
but we do have some information about their number. Important for our purposes is that this 
parametrization clearly distinguishes tempered and nontempered representations. 
These results were inspired by unpublished work of Delorme and Opdam \cite{DeOp2} 
on affine Hecke algebras.

For $\xi = (P,A,\omega ,\chi ) \in \Xi$ we define
\begin{equation}\label{eq:2.9}
\begin{array}{lll@{\qquad}lll}
\nu (\xi ) & = & \log |\chi| & \Sigma (\xi ) & = & \{ \alpha \in
\Sigma (P,A) : \inp{\nu (\xi )}{\alpha}_0 = 0 \} \,, \\
M(\xi ) & = & Z_G (A(\xi ))& A(\xi ) & = & \{ a \in A : \alpha (a)
= 1 \;\forall \alpha \in \Sigma (\xi) \} \,, \\
P(\xi ) & = & P M(\xi ) & \omega (\xi ) & = & I_{P \cap M(\xi
)}^{M (\xi )} \big(\omega \otimes \chi |\chi |^{-1} \big) \,.
\end{array}
\end{equation}
By \cite[Lemme III.2.3]{Wal} $\omega (\xi )$ is a tempered
pre-unitary $M (\xi )$-representation. For $\xi \in \Xi_u$ we simply have
\[
\nu (\xi ) = 0 \,,\, \Sigma (\xi ) = \Sigma (P,A) \,,\, A(\xi ) = A_G \,,\,
P(\xi ) = M(\xi ) = G \text{ and } \omega (\xi ) = I(\xi ) \,.
\]
For general $\xi$ these objects are designed to divide parabolic induction 
into stages, like in \cite[\S XI.9]{KnVo}.
The first stage corresponds to the unitary part of $\chi$ and the
second stage to its absolute value. This is possible since
\begin{align}
I_{P(\xi )}^G \big( |\chi | \otimes \omega (\xi ) \big) &\cong
\mr{Ind}_{P(\xi )}^G \big( \delta_{P(\xi )}^{1/2} \otimes |\chi |
\otimes \omega (\xi ) \big) \nonumber \\
&\cong \mr{Ind}_{P(\xi )}^G \big( \delta_{P(\xi )}^{1/2} \otimes
|\chi | \otimes \mr{Ind}_{P \cap M(\xi )}^{M(\xi )} \big(
\delta_{P \cap M(\xi)}^{1/2} \otimes \chi |\chi|^{-1} \otimes
\omega \big) \big) \nonumber \\
&\cong \mr{Ind}_{P M(\xi )}^G \big( \delta_{P M(\xi )}^{1/2}
\otimes \mr{Ind}_{P \cap M(\xi )}^{M(\xi )} \big( \delta_{P \cap
M(\xi)}^{1/2} \otimes \chi \otimes \omega \big) \big) \label{eq:2.10} \\
&\cong \mr{Ind}_{P M(\xi )}^G \big(  \mr{Ind}_{P \cap M(\xi
)}^{M(\xi )} \big( \delta_{P M(\xi )}^{1/2} \otimes \delta_{P
\cap M(\xi)}^{1/2} \otimes \chi \otimes \omega \big) \big) \nonumber \\
&\cong \mr{Ind}_P^G \big( \delta_P^{1/2} \otimes \chi \otimes
\omega \big) \quad = \quad I(\xi ) \,. \nonumber
\end{align}
Clearly we can transfer the positivity condition from Langlands data to induction data. 
We say that $\xi = (P,A,\omega ,\chi ) \in \Xi^+$ if $(P,A)$ is standard and 
$\log |\chi| \in \bar{\mf a}^{*,+}$. This choice of a ``positive cone" is justified 
by the following result.

\begin{lem}\label{lem:2.10}
Every $\xi \in \Xi$ is $\mc W$-associate to an element of $\Xi^+$.
If $\xi_1 , \xi_2 \in \Xi^+$ are $\mc W$-associate, then the
objects $\Sigma (\xi_i ) ,\, A(\xi_i ) ,\, M(\xi_i ) ,\, P(\xi_i )$ and 
$\nu (\xi_i )$ are the same for $i=1$ and $i=2$, while $\omega (\xi_1 )$ 
and $\omega (\xi_2 )$ are equivalent $M(\xi_i )$-representations.
\end{lem}
\emph{Proof.} As we noted before, every parabolic pair is
conjugate to a standard one. By \cite[Section 1.15]{Hum} every
$W_0$-orbit in $\mf a_0^*$ contains a unique point in a positive
chamber $\mf a^{*,+}$ (for a unique $A \subset A_0$).  This proves
the first claim, and it also shows that
\begin{equation}
\log |\chi_1 | = \log |\chi_2 | \in \mf a_0^* \,.
\end{equation}
Hence the $\nu$'s, $\Sigma$'s, $A$'s, $P$'s and $M$'s are the same
for $i=1$ and $i=2$. If now $w \in \mc W$ is such that $w \omega_1
\cong \omega_2$, then by Theorem \ref{thm:2.4}, applied to
$M(\xi_i )$, there is a unitary intertwiner between $\omega (\xi_1
)$ and $\omega (\xi_2 ). \qquad \Box$
\\[2mm]

The link between these positive induction data and the Langlands 
classification is easily provided:

\begin{prop}\label{prop:2.11}
Let $\xi = (P,A,\omega ,\chi ) \in \Xi^+$.
\begin{description}
\item[a)] Let $\tau$ be an irreducible direct summand of $\omega (\xi )$.
Then $(P(\xi ), A (\xi), \tau ,\nu (\xi )) \in \Lambda^+$.
\item[b)] The irreducible quotients of $I(\xi )$ are precisely the modules
$J (P(\xi ), A (\xi), \tau ,\nu (\xi ))$ with $\tau$ as above, and these
are tempered if and only if $\xi \in \Xi_u$.
\item[c)] The functor $I_{P(\xi )}^G$ induces an isomorphism
$\mr{End}_{M(\xi )} (\omega (\xi )) \cong \mr{End}_G (I(\xi ))$.
\end{description}
\end{prop}
\emph{Proof.}
a) follows directly from the definitions \eqref{eq:2.9}.\\
b) follows from a) and Theorem \ref{thm:2.8}.\\
c) Let $\sigma$ be an irreducible direct summand of $\omega (\xi )$, which is not 
equivalent to $\tau$ as a $M(\xi )$-representation. By Lemma \ref{lem:2.9}.b 
\[
J(P(\xi ), A (\xi), \tau ,\nu (\xi )) \notin 
\mr{JH} \big( I(P(\xi ), A (\xi), \sigma ,\nu (\xi )) \big) \,.
\] 
The proof of \eqref{eq:2.19}, with $I(P(\xi ), A (\xi), \sigma ,\nu (\xi ))$ 
in the role of $M(\lambda )$, shows that 
\[
\mr{Hom}_G \big( I \big( P(\xi ), A (\xi), \tau ,\nu (\xi ) \big) 
,I \big( P(\xi ), A (\xi), \sigma ,\nu (\xi ) \big) \big) = 0 \,.
\]
Now the result follows easily from Lemma \ref{lem:2.9}.c. $\qquad \Box$
\\[2mm]

As announced, we can parametrize Irr$(G)$ with our induction data.

\begin{thm}\label{thm:2.12}
For every $\pi \in \mr{Irr}(G)$ there exists a unique association class\\ 
$\mc W (P,A,\omega ,\chi ) \in \Xi / \mc W$ such that the following equivalent statements hold:
\begin{description}
\item[a)] $\pi$ is equivalent to an irreducible quotient of $I (\xi^+ )$, for some\\ 
$\xi^+ \in \mc W (P,A,\omega ,\chi ) \cap \Xi^+$.
\item[b)] $\pi \in \mr{JH}(I (P,A,\omega ,\chi ))$ and $\norm{\cc_M (\omega )}$ 
is maximal with respect to this property.
\end{description}
\end{thm}
\emph{Proof.} 
a) Let $(P_\pi ,A_\pi ,\sigma ,\nu )$ be the Langlands datum associated to $\pi$. 
Write $\mc W^H, \Xi^H$ etcetera for $\mc W, \Xi$, but now corresponding to 
$H = M_\pi = Z_G (A_\pi )$ instead of $G$. By Theorem \ref{thm:2.5}, applied to $H$, 
there exists a unique association class
\[
\mc W^H \xi^H = \mc W^H (P_F \cap H ,A_F ,\omega^H ,\chi^H ) \in \Xi^H_u / \mc W^H
\]
such that $\sigma$ is a direct summand of 
$I^H (\xi ) = I_{P_F \cap H}^H (\omega^H \otimes \chi^H )$. Put 
\[
\xi^+ = (P_F ,A_F ,\omega^H, \chi^H \cdot \chi_\nu ) \in \Xi^+ \,.
\] 
By Proposition \ref{prop:2.11}.b $\pi \in \mr{JH}(I (\xi^+ )$), and by Lemma \ref{lem:2.10} and 
Theorem \ref{thm:2.5} the class $\mc W \xi^+ \in \Xi / \mc W$ is unique for this property.\\
b) In view of Corollary \ref{cor:2.15} and Lemma \ref{lem:2.10} we may assume that 
$\xi = (P,A,\omega ,\chi ) \in \Xi^+$. Suppose that $\pi$ is not equivalent to a quotient of 
$I(\xi )$, and let $\tau$ be an irreducible summand of $\omega (\xi )$ such that 
$\pi \in \mr{JH}(I (\lambda ))$, with $\lambda = (P(\xi ) ,A(\xi ), \tau ,\nu (\xi )) \in \Lambda^+$. 
By Lemma \ref{lem:2.13} and \eqref{eq:2.21}
\[
\norm{c_{M_F}(\omega^H)} = \norm{\cc_{M_\pi}(\sigma)} <  
\norm{\cc_{M(\xi )}(\tau)} = \norm{\cc_M (\omega)} \,.
\]
This contradiction proves that $\pi$ must be equivalent to a quotient of $I(\xi )$.
Hence $\xi$ satisfies the condition from a), which shows that $\mc W \xi$ is unique.
In particular conditions a) and b) are equivalent. $\qquad \Box$
\\[2mm]

Every standard triple $(P,A,\omega )$ gives a series of finite length $G$-representations 
$I(P,A,\omega ,\chi )$, parametrized by $\chi \in \xnr (M)$. Every such series lies in a 
single Bernstein component. Let Irr$(G)_{(P,A,\omega )}$ denote the set of all 
$\pi \in \mr{Irr}(G)$ which are equivalent to an irreducible quotient of $I(\xi )$, for some 
$\xi \in \Xi^+ \cap \mc W (P,A,\omega, \xnr (M) )$. Theorem \ref{thm:2.12} tells us that 
these subsets form a partition of Irr$(G)$. We note that
\[
\mr{Irr}^t (G)_{(P,A,\omega )} := \mr{Irr}(G)_{(P,A,\omega )} \cap \mr{Irr}^t (G)
\]
consists of the direct summands of the $I(P,A,\omega ,\chi )$ with $\chi \in \xunr (M)$.
As a topological space Irr$(G)_{(P,A,\omega )}$ is usually not Hausdorff, since certain
$\chi \in \xnr {M}$ carry more than one point. The following result shows that the 
singularities depend only on the $\mc W$-action, and can therefore already be detected on 
$\mr{Irr}^t (G)_{(P,A,\omega )}$.

\begin{lem}\label{lem:2.16}
Suppose that $t \to \xi_t = (P,A,\omega ,\chi_t )$ is a path in $\Xi$, and that 
$\mc W_\xi := \{ w \in \mc W : w (\xi_t ) = \xi_t \}$ does not depend on $t$. Then
$| JH (I(\xi_t )) \cap \mr{Irr}(G )_{(P,A,\omega )} |$ does not depend on $t$.
\end{lem}
\emph{Proof.}
For $\chi_t$ moving within $\xunr (M)$ this property can be read 
off from Theorem \ref{thm:2.7}. Take $n \in \mc W_{\xi}$. According  to 
Theorem \ref{thm:2.4} the intertwiner $I(n, \omega \otimes \chi )$ extends
holomorphically to a tubular neighborhood $U$ of $\xunr (M)$ in $\xnr (M)$.
Hence the desired result also holds for paths with $\chi_t \in U \: \forall t$.

Now consider any path as in the statement. For every $t_0$ there is an 
$r \in (-1,0]$ and a neighborhood $T$ of $t_0$ such that 
\[
\chi'_t := \chi_t |\chi_t |^r \in U \: \forall t \in T \,.
\]
By Proposition \ref{prop:2.11}.b
\[
| \mr{JH} (I(\xi_t )) \cap \mr{Irr}(G )_{(P,A,\omega )} | = 
| \mr{JH} (I(\xi'_t )) \cap \mr{Irr}(G )_{(P,A,\omega )} | \,.
\]
We just saw that the right hand side is independent of $t \in T$, hence so is the 
left hand side. $\qquad \Box$ 
\\[2mm]

With these partitions of Irr$ (G)$ and $\mr{Irr}^t (G)$ we can construct filtrations of
corresponding algebras, which will be essential in the next chapter.

\begin{lem}\label{lem:2.14}
Let $\mf s \in \Omega (G)$ be a Bernstein component and pick $K_{\mf s} \in \mr{CO}(G)$ 
as in Proposition \textup{\ref{prop:2.1}.b.} There exist standard triples 
$(P_i ,A_i ,\omega_i )$ and filtrations by two-sided ideals
\begin{align*}
&\mc H (G, K_{\mf s})^{\mf s} = \mc H^{\mf s}_0 \supset \mc H^{\mf s}_1 
\supset \cdots \supset \mc H^{\mf s}_{n_{\mf s}} = 0 \,, \\ 
&\mc S (G, K_{\mf s})^{\mf s} \, = \mc S^{\mf s}_0 \, \supset \, \mc S^{\mf s}_1 \,
\supset \cdots \supset \mc S^{\mf s}_{n_{\mf s}} \, = 0 \,, 
\end{align*}
such that
\begin{itemize}
\item[a.] $\mc H_i^{\mf s} \subset \mc S_i^{\mf s}$ ,
\item[b.] Prim$ \big( \mc H_{i-1}^{\mf s} / \mc H_i^{\mf s} \big) 
\cong \mr{Irr}(G)_{(P_i ,A_i ,\omega_i )}$ ,
\item[c.] Prim$ \big( \mc S_{i-1}^{\mf s} / \mc S_i^{\mf s} \big) \, 
\cong \mr{Irr}^t (G)_{(P_i ,A_i ,\omega_i )}$ .
\end{itemize}
\end{lem}
\emph{Proof.}
By \eqref{eq:2.3} $\mc S (G, K_{\mf s})^{\mf s}$ is a direct summand of 
$\mc S (G, K_{\mf s})$, so Theorem \ref{thm:2.7} assures the existence of finitely many 
(say $n_{\mf s}$) standard triples $(P_i ,A_i ,\omega_i )$ such that
\begin{equation}\label{eq:2.11}
\mc S (G,K_{\mf s})^{\mf s} \cong {\ts \bigoplus_i} \big( C^\infty (\xunr (M_i )) \otimes 
L(\omega_i ,P_i )^{K_{\mf s} \times K_{\mf s}} \big)^{\mc W_i} .
\end{equation}
In view of Theorem \ref{thm:2.15} Prim$ \big( \mc H (G, K_{\mf s})^{\mf s} \big)$ is a union 
of series corresponding to standard triples. Looking at the unitary parts of these series 
and at \eqref{eq:2.11}, we find that
\begin{equation}\label{eq:2.18}
\mr{Prim} \big( \mc H (G, K_{\mf s})^{\mf s} \big) \cong 
{\ts \bigcup_i} \mr{Irr}(G )_{(P_i ,A_i ,\omega_i )} \,.
\end{equation}
Number the triples $(P_i ,A_i ,\omega_i )$ from 1 to $n_{\mf s}$, such that 
\[
\norm{\cc_{M_j} (\omega_j )} \geq \norm{\cc_{M_i} (\omega_i )} \quad \mr{if} \quad j \leq i \,. 
\]
We define
\begin{equation}\label{eq:2.22}
\begin{array}{lll}
\mc H^{\mf s}_i & = & \{ h \in \mc H (G, K_{\mf s})^{\mf s} : \pi (h) = 0 \; \forall \pi \in 
\mr{Irr}(G)_{(P_j ,A_j ,\omega_j )} \,, \forall j \leq i \} \,, \\
\mc S^{\mf s}_i & = & \{ h \in \mc S (G, K_{\mf s})^{\mf s} : \pi (h) = 0 \; \forall \pi \in 
\mr{Irr}^t (G)_{(P_j ,A_j ,\omega_j )} \,, \forall j \leq i \} \,.
\end{array}
\end{equation}
Since the Jacobson closure of $\mr{Irr}^t (G)_{(P_j ,A_j ,\omega_j )}$ in Prim$ (\mc H (G))$ 
contains $\mr{Irr} (G)_{(P_j ,A_j ,\omega_j )} $, we have
$\mc H^{\mf s}_i \subset \mc S^{\mf s}_i$. From \eqref{eq:2.11} we get
\begin{equation}\label{eq:2.23}
 \mc S_{i-1}^{\mf s} / \mc S_i^{\mf s} \cong \big( C^\infty (\xunr (M_i )) \otimes 
L(\omega_i ,P_i )^{K_{\mf s} \times K_{\mf s}} \big)^{\mc W_i} ,
\end{equation}
which shows that c) holds. We claim that
\begin{equation*}
{\ts \bigcup_{j \leq i} } \mr{Irr}(G)_{(P_j ,A_j ,\omega_j )}
\end{equation*}
is closed in the Jacobson topology of $\mc H (G, K_{\mf s})^{\mf s}$.
Its Jacobson closure consists of all irreducible subquotients $\pi$ of 
$I(\xi_j ) = I (P_j , A_j ,\omega_j ,\chi )$, for any $j \leq i$ and $\chi \in \xnr (M_j )$.
Suppose that $\pi \not\in \mr{Irr}(G)_{(P_j ,A_j ,\omega_j )}$. By Theorem \ref{thm:2.12}
$\pi \in \mr{Irr}(G)_{(P ,A ,\omega )}$ for some standard triple with
$\norm{\cc_M (\omega )} > \norm{\cc_{M_j} (\omega_j )}$. In view of \eqref{eq:2.18} we must 
have $w(P,A,\omega) = (P_n ,A_n ,\omega_n )$ for some $w \in \mc W$ and $n \in \mh N$. Then
\[
\norm{\cc_{M_n}(\omega_n )} = \norm{\cc_M (\omega )} > \norm{\cc_{M_j} (\omega_j )} \,,
\]
so $n < j \leq i$, proving our claim. Consequently
\[
\mr{Prim} \big( \mc H_i^{\mf s} \big) \cong 
{\ts \bigcup_{j > i} } \mr{Irr}(G)_{(P_j ,A_j ,\omega_j )} \,,
\]
which easily implies b). $\qquad \Box$ 

\chapter{The noncommutative geometry of reductive $p$-adic groups}

\section{Periodic cyclic homology}
\label{sec:3.1}

The results of Chapters 1 and 2 enable us to compare the $K$-theory and
the periodic cyclic homology of reductive $p$-adic groups. Before 
proceeding we recall that $HP_*$ is continuous in certain situations.
Here and in the next results we will freely use the notations from Chapter 2.

\begin{thm}\label{thm:3.1}
\begin{description}
\item[a)] Suppose that $A = \lim_{i \to \infty} A_i$ is an inductive limit of algebras
and that there exists $N \in \mh N$ such that $HH_n (A_i ) = 0 \, \forall n > N , 
\forall i$. Then $HP_* (A) \cong \lim_{i \to \infty} HP_* (A_i )$.
\item[b)] Suppose that  $B = \lim_{i \to \infty} B_i$ is a strict inductive limit of nuclear
Fr\'echet algebras and that there exists $N \in \mh N$ such that 
$HH_n (B_i ,\hot ) = 0 \, \forall n > N , \forall i$. 
Then $HP_* (B ,\hot_{\mh C} ) \cong \lim_{i \to \infty} HP_* (B_i ,\hot)$.
\end{description}
\end{thm}
\emph{Proof.}
According to \cite[Theorem 1.93]{Mey3} the completed bornological tensor product
agrees with Grothendieck's completed inductive tensor product for strict inductive limits 
of nuclear Fr\'echet spaces. This identifies b) with \cite[Theorem 3]{BrPl1}.
Part a) is just the simpler algebraic version of this result, which can also be found in
\cite[Proposition 2.2]{Nis} $\qquad \Box$
\\[2mm]
\begin{thm}\label{thm:3.2}
Let $\mf s \in \Omega (G)$ be a Bernstein component and let $K_{\mf s} \in \mr{CO}(G)$
be as in Proposition \textup{\ref{prop:2.1}.b.}
\begin{description}
\item[a)] The Chern character for $\mc S (G,K_{\mf s} )^{\mf s}$ induces an isomorphism
\[
K_* (C_r (G)^{\mf s} ) \otimes_{\mh Z} \mh C \to HP_* (\mc S (G)^{\mf s} ,\hot_{\mh C} ) \,.
\]
\item[b)] The direct sum of these maps, over all $\mf s \in \Omega (G)$, 
is a natural isomorphism
\[
K_* (C_r (G)) \otimes_{\mh Z} \mh C \to HP_* (\mc S (G),\hot_{\mh C} ) \,.
\]
\end{description}
\end{thm}
\emph{Proof.}
By \eqref{eq:2.11} and Theorem \ref{thm:1.6} the Chern character
\begin{equation}\label{eq:3.24}
ch \otimes \mr{id}_{\mh C} : K_* \big( \mc S (G,K_{\mf s} )^{\mf s} \big) \otimes_{\mh Z} \mh C
\to HP_* \big( \mc S (G,K_{\mf s} )^{\mf s} \big)
\end{equation}
is an isomorphism. Furthermore there are natural isomorphisms
\begin{equation}\label{eq:3.1}
\begin{aligned}
K_* (C_r^* (G)) \, &\cong K_* \Big( \varinjlim_{\mf S} 
\bigoplus_{\mf s \in \mf S} C_r^* (G)^{\mf s} \Big) \\
& \cong \bigoplus_{\mf s \in \Omega (G)} K_* \left( C_r^* (G)^{\mf s} \right) \\
& \cong \bigoplus_{\mf s \in \Omega (G)} \varinjlim_{K \in \mr{CO} (G)} 
K_* \left( C_r^* (G,K)^{\mf s} \right) \\
& \cong \bigoplus_{\mf s \in \Omega (G)} 
K_* \big( C_r^* (G,K_{\mf s} )^{\mf s} \big)\\
&\cong \bigoplus_{\mf s \in \Omega (G)} 
K_* \big( \mc S (G,K_{\mf s} )^{\mf s} \big) \,.
\end{aligned}
\end{equation}
Here we used respectively \eqref{eq:2.13}, \eqref{eq:2.14},
Proposition \ref{prop:2.1}.d and that $K_*$ is invariant for passing 
to holomorphically closed dense Fr\'echet subalgebras. Similarly,
\begin{equation}\label{eq:3.2}
\begin{aligned}
HP_* (\mc S (G), \hot_{\mh C} ) \: & \cong 
\varinjlim_{K \in \mr{CO} (G)} HP_* (\mc S (G,K), \hot )\\
&\cong \bigoplus_{\mf s \in \Omega (G)} \varinjlim_{K \in \mr{CO} (G)}
HP_* (\mc S (G,K)^{\mf s}, \hot )\\
& \cong \bigoplus_{\mf s \in \Omega (G)} 
HP_* (\mc S (G)^{\mf s}, \hot_{\mh C} )\\
&\cong \bigoplus_{\mf s \in \Omega (G)} 
HP_* (\mc S (G, K_{\mf s})^{\mf s}, \hot ) \,.
\end{aligned}
\end{equation}
For the first and third isomorphisms we need Theorem \ref{thm:3.1}.b,
which we may indeed apply due to Corollary \ref{cor:2.3}.b. The second and
fourth isomorphisms rely on Proposition \ref{prop:2.1}.c. 

Finally we combine the last lines of  \eqref{eq:3.1} and \eqref{eq:3.2} with \eqref{eq:3.24}.
$\qquad \Box$ \\[2mm]

Altogether this theorem involves quite a few steps, but one is always
guided by the general principle that algebras with the same spectrum should 
have closely related invariants. 

Next we will prove the comparison theorem for the periodic cyclic homology
of reductive $p$-adic groups. This result was suggested in 
\cite[Conjecture 8.9]{BHP2} and in \cite[Section 4]{ABP}. It is more difficult than 
Theorem \ref{thm:3.2}, precisely because the above principle does not apply.
We remark that Nistor \cite[Theorem 4.2]{Nis} already gave a rather explicit 
description of $HP_* (\mc H (G))$.

\begin{thm}\label{thm:3.3}
The inclusions $\mc H (G)^{\mf s} \to \mc S (G)^{\mf s}$ for $\mf s \in \Omega (G)$
induce isomorphisms
\begin{align*}
HP_* ( \mc H (G)^{\mf s} ) &\isom HP_* ( \mc S (G)^{\mf s}, \hot_{\mh C} ) \,, \\
HP_* ( \mc H (G) ) &\isom HP_* ( \mc S (G),\hot_{\mh C} ) \,.
\end{align*}
\end{thm}
\emph{Proof.}
It follows from Theorem \ref{thm:3.1}.a, Corollary \ref{cor:2.3}.a and Proposition 
\ref{prop:2.1}.b that 
\begin{equation}\label{eq:3.4}
\begin{aligned}
HP_* (\mc H (G) ) \: &\cong 
\bigoplus_{\mf s \in \Omega (G)} 
HP_* (\mc H (G)^{\mf s} )\\
&\cong \bigoplus_{\mf s \in \Omega (G)} \varinjlim_{K \in \mr{CO}(G)}
HP_* (\mc H (G,K)^{\mf s} )\\
&\cong \bigoplus_{\mf s \in \Omega (G)} 
HP_* (\mc H (G, K_{\mf s})^{\mf s} )\\
\end{aligned}
\end{equation}
Together with \eqref{eq:3.2} this reduces the proof to showing that for every
Bernstein component $\mf s \in \Omega (G)$ the inclusion
\begin{equation}
\mc H (G,K_{\mf s} )^{\mf s} \to \mc S (G,K_{\mf s} )^{\mf s} 
\end{equation}
induces an isomorphism on \pch \!. Applying Lemma \ref{lem:1.3} to the filtrations
from Lemma \ref{lem:2.14}, this will follow from the next result.

\begin{prop}\label{prop:3.8}
Let $\mc H_i^{\mf s}$ and $\mc S_i^{\mf s}$ be as in \eqref{eq:2.22}. The inclusion map
\[
\mc H_{i-1}^{\mf s} / \mc H_i^{\mf s} \to \mc S_{i-1}^{\mf s} / \mc S_i^{\mf s}
\]
induces an isomorphism on \pch \!.
\end{prop}
\emph{Proof.}
We would like to copy the proof of Theorem \ref{thm:1.7} with $\Gamma = \mc W_i$ and
\[
\begin{array}{lll@{,\qquad}lll}
A\al & = & \mc H^{\mf s}_{i-1} / \mc H^{\mf s}_i & X & = & \xnr  (M_i ) \,, \\
A\sm & = & \mc S^{\mf s}_{i-1} / \mc S^{\mf s}_i & X' & = & \xunr  (M_i ) \,.
\end{array}
\]
By Lemmas \ref{lem:2.16} and \ref{lem:2.14} the primitive ideal spectra fit into the framework
of Section \ref{sec:1.4}. On the smooth side everything works fine, but on the algebraic
side we have to take into account that in general
\begin{itemize}
\item the intertwiners $u_w \; (w \in \mc W_i )$ are only rational on $\xnr (M_i )$,
\item the center of $\mc H^{\mf s}_{i-1} / \mc H^{\mf s}_i$ is ``too small", in the sense that
it does not contain all polynomial functions on $\xnr (M_i ) / \mc W_i$.
\end{itemize}
In fact every element of 
\[
Z \big( \mc H^{\mf s}_{i-1} / \mc H^{\mf s}_i \big) \subset \mc O (\xnr (M_i ) / \mc W_i )
\]
necessarily vanishes on the singularities of the intertwiners, because these points of 
$\xnr (M_i )$ carry representations from Irr$ (G )_{(P_j ,A_j ,\omega_j )}$ with $j < i$.

Let $\mc L$ be the collection of all the irreducible components of all the $X^H$, with $H$ 
running over all subsets of $\mc W_i$. We note that these are all of the form $\chi T$, 
with $\chi \in \xunr (M_i )$ and $T$ an algebraic subtorus of $\xnr (M_i )$. 
We write $T_u = T \cap \xunr (M_i)$. Let $\mc L_p$ be the subset of $\mc L$ consisting
of elements of dimension $\leq p$, and define the $\mc W_i$-stable subvarieties
\[
\begin{array}{lll}
X_p & := & \bigcup_{\chi T \in \mc L_p} \, \chi T \,, \\
X'_p & := & X_p \cap \xunr (M_i ) = 
\bigcup_{\chi T \in \mc L_p} \, \chi T_u \,.
\end{array}
\]
Let $V$ be the finite dimensional vector space $I^G_{P_i} (E_i )^{K_{\mf s}}$ and
consider the ideals 
\[
\begin{array}{lll}
I_p & := & \{ h \in \mc H^{\mf s}_{i-1} / \mc H^{\mf s}_i : 
I(P_i ,A_i ,\omega_i ,\chi )(h) = 0 \; \forall \chi \in X_p \} \,, \\
J_p & := & \{ h \in \mc S^{\mf s}_{i-1} / \mc S^{\mf s}_i : 
I(P_i ,A_i ,\omega_i ,\chi )(h) = 0 \; \forall \chi \in X'_p \} \,.
\end{array}
\]
It follows from \eqref{eq:2.23} that 
\[
\begin{array}{lll}
J_{p-1} / J_p & = & C_0^\infty (X'_p , X'_{p-1} ; \mr{End}(V) )^{\mc W_i} \,, \\
I_{p-1} / I_p & \subset & \mc O_0 (X_p , X_{p-1} ; \mr{End}(V) )^{\mc W_i} \,.
\end{array}
\]
The Jacobson topology makes Prim$ (A_{alg} / I_p )$ into a complex scheme,
which need not be separable. Let 
\[
\theta_p : \mr{Prim}(A_{alg} / I_p ) \to X_p / \mc W_i
\]
be the natural morphism. Take $C \in \mc L_p$ and consider
$Y_C = \theta_p^{-1} ((\mc W_i C) / \mc W_i)$. From Lemma \ref{lem:2.16} we see that 
$Y_C$ is a disjoint union of copies of $(\mc W_i C) / \mc W_i $, of which certain subvarieties 
are given a higher multiplicity. Moreover, the only thing that can cause inseparability on
Prim$ (A\al / I_p )$ is a jump in the isotropy group of $\mc W_i$ acting on $C$. By construction
of $X_p$ this implies that all inseparable points lie in $\theta_p^{-1} (X_{p-1} / \mc W_i )$, 
so the corresponding primitive ideals contain $I_{p-1} / I_p$. 

We define $Y_p$ as the (maximal) separable quotient of Prim$ (A\al / I_p )$. From the above 
description we see that it is a complex algebraic variety on which $\theta_p$ is well-defined. 
Let $Z_p$ be its closed subvariety corresponding to $X_{p-1} / \mc W_i$. Then
\[
\mr{Prim}(I_{p-1} /I_p ) = Y_p \setminus Z_p \,,
\]
and this would agree with \eqref{eq:1.20} if all the intertwiners were polynomial. We have
\[
\begin{array}{lllll}
Y'_p & := & Y_p \cap \theta_p^{-1} (X'_p / \mc W_i ) & = & 
\text{Hausdorff quotient of Prim} (A\sm / J_p ) \,, \\
Z'_p & := & Z_p \cap \theta_p^{-1} (X'_p / \mc W_i ) & = & 
\{ J/ \! \sim \, \in Y'_p : Z(J_{p-1} / J_p ) \subset J \} \,, \\
Y'_p \setminus Z'_p & = & \mr{Prim}(J_{p-1} /J_p ) & = & \mr{Prim} \big( Z(J_{p-1} /J_p ) \big) \,.
\end{array}
\]
With these notations the singularities of the $u_w$ form a closed subvariety of $Y_p$, 
disjoint from $Y'_p$. From page \pageref{eq:1.21} we know that there are natural isomorphisms
\begin{equation}
HP_* (J_{p-1} / J_p ) \leftarrow HP_* \big( Z(J_{p-1} / J_p ) \big) =
HP_* \big( \tilde C_0^\infty (Y'_p ,Z'_p) \big) \cong H^{[*]} (Y'_p ,Z'_p ) \,.
\end{equation}
For our comparison we will construct something similar on the algebraic side.

\begin{lem}\label{lem:3.9}
There exists a natural map
\[
HP_* (I_{p-1} /I_p ) \to HP_* (\mc O_0 (Y_p ,Z_p )) \,. 
\]
\end{lem}
\emph{Proof.}
The main problem is the absence of a natural algebra homomorphism between $I_{p-1} / I_p$ and 
$\mc O_0 (Y_p ,Z_p)$. As a substitute we will use the generalized
trace map on the periodic cyclic bicomplex.

Let $C$ be an irreducible component of $X_p$, which is not contained in $X_{p-1}$. Pick
any point $\chi \in C \setminus (C \cap X_{p-1})$ and define
\[
\mc W_C := \{ w \in \mc W : w \chi = \chi \} \,.
\]
This depends only on $C$, not on the choice of $\chi$. Since $C$ is of the form ``algebraic
subtorus translated by a unitary element", $C \cap \xunr (M_i )$ is a (nonempty) real form 
of $C$. According to Theorem \ref{thm:2.4} the intertwiners $I (w,\omega \otimes \chi )$ with 
$w \in \mc W_C$ are regular on $C \cap \xunr (M_i )$, so they are regular on a nonempty 
Zariski-open subset of $C$. But $I(w ,\omega \otimes \chi ) \in GL ( \mr{End}_{\mh C}(V))$ 
has finite order for all $\chi \in C$, so it is in fact regular on the whole of $C$. In particular, 
for every $\chi \in C$ there is a canonical decomposition of $V$ into isotypical projective 
$\mc W_C$-representations:
\begin{equation}\label{eq:3.6}
I(P_i ,A_i ,\omega_i ,\chi )^{K_{\mf s}} = V = V_1^\chi \oplus \cdots \oplus V_{n_C}^\chi \,.
\end{equation}
Since $\chi \mapsto I(w,\omega \otimes \chi)$ is polynomial on $C$, the type of $V$ as a 
projective $\mc W_C$-representation is independent of $\chi \in C$. Moreover the corresponding 
projections $p(C,n,\chi) \in \mr{End}_{\mh C} (V)$ are polynomial in $\chi$.

For $\chi \in C \cap (X'_p \setminus X'_{p-1} )$ Theorem \ref{thm:2.7} assures that
\[
\mr{End}_{\mc W_C} \big( I(P_i ,A_i ,\omega_i ,\chi )^{K_{\mf s}} \big) = 
\{ I(P_i ,A_i ,\omega_i ,\chi) (h) : h \in \mc S^{\mf s}_{i-1} / \mc S^{\mf s}_i \} \,,
\]
so the summands of \eqref{eq:3.6} are in bijection with the irreducible constituents of the
$\mc H^{\mf s}_{i-1} / \mc H^{\mf s}_i$-module $I(P_i ,A_i ,\omega_i ,\chi )^{K_{\mf s}}$.
Together with Lemma \ref{lem:2.16} this shows that for every 
$\chi \in C \cap (X_p \setminus X_{p-1})$ and every 
$\pi \in \mr{JH}_{A\al} \big( I(P_i ,A_i ,\omega_i ,\chi )^{K_{\mf s}} \big)$ there is
a unique direct summand of \eqref{eq:3.6} in which $\pi$ appears. We remark that  
$\mr{JH}_G ( I(P_i ,A_i ,\omega_i ,\chi ) )$ may be larger, but the 
additional elements do not belong to Irr$ (G )_{(P_i ,A_i ,\omega_i )}$.

For any $w \in \mc W_i$ and any $n \leq n_C$ the intertwiner $I(w,\omega \otimes \chi)$ maps
$\mr{End}_{\mh C}(V_n^\chi )$ to $\mr{End}_{\mh C}(V_{n'}^{w \chi} )$, for some isotypical
component $V_{n'}^{w \chi}$ of the projective $(w \mc W_C w^{-1})$-representation 
$I(P_i ,A_i ,\omega_i ,w \chi )^{K_{\mf s}}$. Write 
\[
\tilde Y_p := {\ts \bigsqcup_C} \{ 1,\ldots ,n_C \} \times C \quad , \quad 
\tilde Z_p := {\ts \bigsqcup_C} \{ 1,\ldots ,n_C \} \times (C \cap X_{p-1}) \,,
\]
and let $\mc W_i$ act on these spaces by sending $(n,C, \chi)$ to $(n',w C ,w \chi)$.
The above amounts to a bijection
\begin{equation}\label{eq:3.3}
\mr{Prim}(I_{p-1}/I_p ) = Y_p \setminus Z_p \to (\tilde Y_p \setminus \tilde Z_p ) / \mc W_i \,.
\end{equation}
Furthermore we can make $\mc W_i$ act on the algebra 
\[
B_p := \mc O_0 (\tilde Y_p ,\tilde Z_p ) \otimes \mr{End}(V)
\]
by the natural combination of the actions on $I_{p-1}/I_p$ and on $\tilde Y_p$. There is a
morphism of finite type $Z \big( \mc H (G,K_{\mf s})^{\mf s} \big)$-algebras
\begin{align*}
& \phi : I_{p-1}/I_p \to B_p^{\mc W_i} \,, \\
& \phi (h) (C,n,\chi) = p(C,n,\chi ) \, I(P_i ,A_i ,\omega_i ,\chi )(h) \, p(C,n,\chi )\,.
\end{align*}
By \eqref{eq:3.3} the central subalgebra
\[
Z(B_p )^{\mc W_i} = \mc O_0 (\tilde Y_p ,\tilde Z_p )^{\mc W_i} \subset B_p^{\mc W_i}
\]
satisfies
\[
\mr{Prim} \big( Z(B_p )^{\mc W_i} \big) = (\tilde Y_p \setminus \tilde Z_p ) / \mc W_i
\cong \mr{Prim}(I_{p-1}/I_p ) \,.
\]
Now we are in the right position to apply the generalized trace map. 
Recall \cite[Section 1.2]{Lod} that, for any algebra $B$, this is a collection of linear maps 
\[
\mr{tr}_m : (B \otimes \mr{End} (V) )^{\otimes m} \to B^{\otimes m}
\]
which together form chain maps on the standard complexes computing Hochschild
and (periodic) cyclic homology. The tracial property can be formulated as
\begin{equation}\label{eq:3.5}
\mr{tr}_n \circ \mr{Ad}(b)^{\otimes n} = \mr{tr}_n
\end{equation}
for all invertible elements $b$ in the multiplier algebra of $B$. The induced maps on 
homology are natural and inverse to the maps induced by the inclusion
\begin{equation}\label{eq:3.11}
B \to B \otimes \mr{End}(V) \;:\; b \to e b e
\end{equation}
with $e \in \mr{End}(V)$ an idempotent of rank one. Consider the composition 
\[
\mr{tr}_m \circ \phi^{\otimes m} : (I_{p-1}/I_p )^{\otimes m} \to Z(B_p )^{\otimes m} \,.
\]
In view of \eqref{eq:3.5} and \eqref{eq:2.7} the image is contained in 
$\big( ( Z(B_p ) )^{\otimes m} \big)^{\mc W_i}$. The periodic cyclic bicomplex 
\cite[Section 5.1]{Lod} with terms $\big( ( Z(B_p ) )^{\otimes m} \big)^{\mc W_i}$ 
has homology $HP_* (Z(B_p ))^{\mc W_i}$, because $\mc W_i$ is finite and 
acts by algebra automorphisms. Thus we obtain a natural map
\begin{equation}\label{eq:3.25}
\tau := HP_* (\mr{tr} \circ \phi ) : HP_* (I_{p-1}/I_p ) \to HP_* (Z(B_p ))^{\mc W_i} \,.
\end{equation}
By Theorem \ref{thm:1.4}.a there are natural isomorphisms
\begin{equation}\label{eq:3.28}
\begin{split}
HP_* (Z(B_p ))^{\mc W_i} & \leftarrow HP_* \big( Z(B_p )^{\mc W_i} \big) = 
HP_* \big(  \mc O_0 (\tilde Y_p ,\tilde Z_p )^{\mc W_i} \big) \\
& = HP_* (\mc O_0 (Y_p ,Z_p )) \cong H^{[*]}(Y_p ,Z_p) \,.
\end{split}
\end{equation}
The combination of \eqref{eq:3.25} and \eqref{eq:3.28} yields the required map. 
$\qquad \Box$ 
\\[2mm]
\begin{lem}\label{lem:3.10}
The map 
\[
HP_* (I_{p-1} /I_p ) \to HP_* (\mc O_0 (Y_p ,Z_p ))
\]
constructed in Lemma \textup{\ref{lem:3.9}} is an isomorphism.
\end{lem}
\emph{Proof.} We will prove the equivalent statement that \eqref{eq:3.25} is bijective.

For $q \in \mh N$, let $Y_{p,q} \subset \mr{Prim}(I_{p-1}/I_p)$ be the collection of all primitive
ideals for which the corresponding irreducible $I_{p-1}/I_p $-module has dimension $\leq q$. 
According to \cite[p. 328]{KNS} these are Jacobson-closed subsets, and they give rise to 
the standard filtration
\begin{equation}\label{eq:3.26}
\begin{split}
& I_{p-1}/I_p = L_0 \supset L_1 \supset \cdots \supset L_{\dim V} \,, \\
& L_q := {\ts \bigcap_{J \in Y_{p,q}} } J \,. \\
\end{split}
\end{equation}
Since every irreducible $I_{p-1}/I_p$-module is a subquotient of $V \,, L_{\dim V}$ is the 
Jacobson radical of $I_{p-1}/I_p$. We would like to apply the results of \cite{KNS} to 
\eqref{eq:3.26}, but this is not directly possible because $I_{p-1}/I_p$ does not have a unit. 
Therefore we consider the unital finite type algebra
\[
\begin{array}{lll}
A & := & Z \big( \mc H (G,K_{\mf s})^{\mf s} \big) \oplus M_2 (I_{p-1}/I_p ) \,, \\
(f \oplus h) (f' \oplus h') & := & (ff' \oplus fh' + f' h + hh' ) \,.
\end{array}
\]
Since $\mc H_i^{\mf s}$ and $\mc H_{i-1}^{\mf s}$ are ideals in $\mc H (G,K_{\mf s})^{\mf s}$ and
every $I_p \subset \mc H^{\mf s}_{i-1} / \mc H^{\mf s}_i$ is defined in terms of $\mc W_i$-stable
subvarieties of $\xnr (M_i )$, this multiplication is well-defined.
The standard filtration of $A$ is
\[
A \supset M_2 (L_0) \supset M_2 (L_0 ) \supset M_2 (L_1) \supset M_2 (L_1) \supset \cdots \supset 
M_2 (L_{\dim V}) \,.
\]
According to \cite[Proposition 1 and Theorem 10]{KNS}, applied to the algebra $A$ and its ideal
$M_2 (L_0 ) = M_2 (I_{p-1}/I_p )$, there are natural isomorphisms
\begin{equation}\label{eq:3.27}
HP_* (M_2 (L_{q-1}) / M_2 (L_q) ) \cong H^{[*]} (Y_{p,q} \cup Z_p ,Y_{p,q-1} \cup Z_p ) \,.
\end{equation}
Let $\tilde Y_{p,q}$ be the inverse image of $Y_{p,q}$ under the map $\tilde Y_p \setminus \tilde Z_p
\to Y_p \setminus Z_p$ from \eqref{eq:3.3}. To \eqref{eq:3.26} corresponds the filtration
\begin{align*}
& Z(B_p) = F_0 \supset F_1 \supset \cdots \supset F_{\dim V} = 0 \,, \\
& F_q := \mc O_0 (\tilde Y_p , \tilde Y_{p,q} \cup Z_p ) \,.
\end{align*}
By construction \eqref{eq:3.25} induces a map
\[
\tau_q : HP_* (L_{q-1}/L_q ) \to HP_* (F_{q-1}/F_q )^{\mc W_i} \,,
\]
which is natural with respect to morphisms of the underlying varieties.
By Theorem \ref{thm:1.4}.a the right hand side is isomorphic to
\begin{equation}\label{eq:3.29}
\begin{split} HP_* (\mc O_0 (\tilde Y_{p,q} \cup \tilde Z_p 
,\tilde Y_{p,q-1} \cup \tilde Z_p ) )^{\mc W_i} 
& \cong HP_* (\mc O_0 (Y_{p,q} \cup Z_p ,Y_{p,q-1} \cup Z_p )) \\
& \cong H^{[*]} (Y_{p,q} \cup Z_p ,Y_{p,q-1} \cup Z_p ) \,.
\end{split}
\end{equation}
Since $\tau_q$ is natural, \eqref{eq:3.27} and \eqref{eq:3.29} show that it is an isomorphism.
Although $\tau$ and $\tau_q$ are not induced by algebra homomorphisms, they do come from maps of
the appropriate periodic cyclic bicomplexes, which is enough to ensure that they are compatible
with the connecting maps from \eqref{eq:1.19}. Hence we can use a variation on Lemma \ref{lem:1.3}, 
where the role of $HP_* (\phi )$ is played by $\tau_*$.
This leads to the conclusion that \eqref{eq:3.25} is an isomorphism. $\qquad \Box$
\\[2mm]

We return to the proof of Proposition \ref{prop:3.8}.
Lemmas \ref{lem:3.9} and \ref{lem:3.10} allow us to write down the following diagram:
\[
\begin{array}{ccccc}
HP_* (I_{p-1}/I_p ) & \rightarrow & HP_* \big( \mc O_0 (Y_p ,Z_p ) \big) 
& \cong & H^{[*]} (Y_p , Z_p ) \\
\downarrow {\scs (1)} & & \downarrow {\scs (2)} & & \downarrow {\scs (3)} \\
HP_* (J_{p-1}/J_p ) & \leftarrow & HP_* \big( \tilde C_0^\infty (Y'_p , Z'_p ) \big) 
& \cong &H^{[*]} (Y'_p ,Z'_p ) 
\end{array}
\]
By Lemma \ref{lem:2.16} $(Y'_p ,Z'_p )$ is a deformation retract of $(Y_p ,Z_p )$, so (3) and (2)
are isomorphisms. Unfortunately the diagram does not commute, because unlike \eqref{eq:3.11} 
$\tilde C_0^\infty (Y'_p , Z'_p ) \to J_{p-1} / J_p$ is not a ``rank one" inclusion. However, 
all the vector spaces in the diagram decompose as direct sums over the components of 
$Y_p \setminus Z_p$, which were labelled $(C,n)$ on page \pageref{eq:3.3}. For every such 
component the diagram commutes up to a scalar factor, namely the dimension of the corresponding 
module $V^x_n$ from \eqref{eq:3.6}. Therefore the diagram does show that the arrow (1) is an 
isomorphism. This concludes the proofs of Proposition \ref{prop:3.8} and Theorem \ref{thm:3.3}.
$\qquad \Box$ \\[2mm]

We remark that a slightly simpler version of the above proof also works for 
affine Hecke algebras and their Schwartz completions, see \cite{Sol2}. 
\vspace{4mm}

\section{Example: $SL_2 (\mh Q_p )$}
\label{sec:3.2}

To clarify the proof of Theorem \ref{thm:3.3} we show in some detail what it
involves in the simplest case. This section is partly based on the related 
calculations in \cite{BHP1} and \cite[Section 6.1]{Sol3}.

Let $p$ be an odd prime, $\mh Q_p$ the field of
$p$-adic numbers, $\mh Z_p$ the ring of $p$-adic integers and $p \mh Z_p$
its unique maximal ideal. We consider the reductive group
$G = SL_2 (\mh Q_p )$ with the maximal torus
$A = \big\{ \matje{a}{0}{0}{a^{-1}} \,: a \in \mh Q^\times_p \big\}$
and the minimal parabolic subgroup
$P = \big\{ \matje{a}{b}{0}{a^{-1}} \,: a \in \mh Q^\times_p\,, b \in \mh Q_p \big\} $.
We have 
\begin{align*}
& M = Z_G (A)  = A \,, \\
& W = N_G (A) / Z_G (A) = \big( A \cup \matje{0}{-1}{1}{0} A \big) / A \,.
\end{align*}
The Iwahori subgroup is
$K = \big\{ \matje{a}{b}{c}{d} \,: a,b,d \in \mh Z_p \,, c \in p \mh Z_p \big\} $.
In this situation $\mc H (G,K)$ is Morita equivalent to $\mc H (G )^{\mf s}$,
where $\mf s = [M,\sigma ] \in \Omega (G)$ is the Borel component, 
corresponding to the trivial representation $(\sigma ,E)$ of $M$. According to
\cite{IwMa} $\mc H (G,K)$ is isomorphic to the Iwahori--Hecke algebra 
$\mc H (A_1 ,p)$ of type $A_1$ with parameter $p$. Furthermore $\mc S (G,K)$
is isomorphic to the Schwartz completion $\mc S (A_1 ,p)$ of $\mc H (A_1 ,p)$,
see \cite{DeOp1}.

We identify $\xnr (M)$ with $\mh C^\times$ by evaluation at $\matje{p}{0}{0}{p^{-1}}$ 
For almost all $\chi \in \xnr (M)$ the $G$-representation $I(P,A,\sigma ,\chi )$ 
is irreducible, so the separated quotient of Prim$( \mc H (G )^{\mf s}$ is 
\[
\xnr (M) / W \cong \mh C^\times / (z \sim z^{-1}) \,.
\]
The $K$-invariant part $I' (\chi)$ of $I(P,A,\sigma ,\chi )$ is a two-dimensional 
$\mc H (G,K)$-module with underlying vector space $V_\chi = V := I_P^G (E)^K$.
The intertwining operator $u_w (\chi )$ has rank one if $\chi^2 = p^{\pm 1}$
and is invertible for all other $\chi \in \xnr (M)$. 
More precisely the homomorphism
\[
I' (\chi ) : \mc H (G,K) \to \mr{End}_{\mh C} (V_\chi )
\]
is surjective for the generic $\chi$ and has image conjugate to
$\{ \matje{a}{b}{0}{d} \,: a,b,d \in \mh C \}$ for special $\chi$. Therefore Prim$
(\mc H (G,K))$ has only two pairs of nonseparated points, at
\begin{equation}
W \chi = p^{\pm 1/2} \quad \text{and} \quad W \chi = -p^{\pm 1/2} \,.
\end{equation}
The $\mc H (G,K)$-modules $V_\chi$ with $\chi \in \xunr (M)$
extend continuously to $\mc S (G,K)$-modules. Besides that, $\mc S (G,K)$
admits precisely two inequivalent one-dimensional square-integrable modules,
namely the irreducible submodules of $V_{p^{1/2}}$ and of $V_{-p^{1/2}}$. Hence 
Prim$ (\mc S(G,K))$ consists of two isolated points (say $\delta_+$ and 
$\delta_-$) and a copy of $\xunr (M) / W \cong [-1,1]$.
In the filtrations 
\begin{align*}
&\mc H (G, K_{\mf s})^{\mf s} = \mc H^{\mf s}_0 \supset \mc H^{\mf s}_1 
\supset \mc H^{\mf s}_2 = 0 \,, \\ 
&\mc S (G, K_{\mf s})^{\mf s} \, = \,\mc S^{\mf s}_0 \supset \, \mc S^{\mf s}_1  
\supset \, \mc S^{\mf s}_2 = 0 \,,
\end{align*}
we have
\begin{align*}
& \mc S_0^{\mf s} \cong C^\infty (\xnr (M))^W \otimes 
\mr{End}_{\mh C} (V) \; \oplus \; \mh C \; \oplus \; \mh C \,,\\
& \mc S_1^{\mf s} \cong C^\infty (\xnr (M))^W \otimes 
\mr{End}_{\mh C} (V) \,,\\
& \mc H_1^{\mf s} = \ker (\delta_+ ) \cap \ker (\delta_- )  \subset
\mc O \big( \xnr (M) ;\mr{End}_{\mh C} (V) \big)^W \,, \\
& \mc H_0^{\mf s} / \mc H_1^{\mf s} \cong \mc S_0^{\mf s} / \mc S_1^{\mf s}
\cong \mr{End}_{\mh C}(\delta_+ ) \oplus \mr{End}_{\mh C}(\delta_- ) \cong  
\mh C \oplus \mh C \,.
\end{align*}
The tricky step is to see that $HP_* (\mc H_1^{\mf s}) \cong 
HP_* (\mc S_1^{\mf s})$. Clearly
\begin{align*}
& \mr{Prim}(\mc S_1^{\mf s}) \cong \xunr (M) /W \,, \\
& Z (\mc S_1^{\mf s}) \cong  C^\infty (\xunr (M))^W \,,\\
& HP_* (\mc S_1^{\mf s}) \cong HP_* (Z(\mc S_1^{\mf s})) \cong
H^{[*]}( \xunr (M) / W) \cong \check H^* ([-1,1] ; \mh C) \,.
\end{align*}
However the image of 
\[
I' \big( \pm p^{1/2} \big) : \mc H_1^{\mf s} \to \mr{End}_{\mh C} ( V_{\pm p^{1/2}} )
\]
is not $M_2 (\mh C )$, but it is conjugate to 
$\{ \matje{0}{b}{0}{d} \,: b,d \in \mh C \}$. Therefore 
\[
Z (\mc H_1^{\mf s}) \cong \mc O_0 \big(\xnr (M) / W , 
\{ p^{\pm 1/2} ,-p^{\pm 1/2} \} \big) \,,
\]
even though Prim$ (\mc H_1^{\mf s}) \cong \xnr (M) / W$. 
Consider the diagram
\[
\begin{array}{ccccc}
\ker I' \big( p^{1/2} \big) \cap \ker I' \big( -p^{1/2} \big) & \to & \mc H_1^{\mf s} & 
\to & \mr{End} ( V_{p^{1/2}}) \oplus \mr{End} (V_{-p^{1/2}} ) \\
\uparrow & & \downarrow \mr{tr} & & \uparrow \\
\! \mc O_0 \big(\xnr (M) / W , W \{ \pm p^{1/2} \} \big) & \to &
\mc O \big(\xnr (M) / W \big) & \to & \mc O \big( \{ p^{1/2} ,-p^{1/2} \} \big) 
\end{array}
\]
The upward arrows identify the centers of the respective algebras. These
morphisms are spectrum preserving, so they induce isomorphisms on
periodic cyclic homology. The downward arrow is the (generalized) trace map, 
which induces a map
\[
HP_* (\mr{tr}) : HP_* (\mc H_1^{\mf s}) \to HP_* (\mc O (\xnr (M) ) )^W
\cong HP_* (\mc O (\xnr (M)/W)) \,.
\]
Now we apply the functor $HP_*$ to the entire diagram, and we replace the
downward arrow by $\frac{1}{2} HP_* (\mr {tr})$. The resulting diagram 
commutes and shows that $HP_* (\mr{tr})$ is a natural isomorphism.
From the commutative diagram 
\[
\begin{array}{ccccc}
HP_* (\mc H_1^{\mf s}) & \xrightarrow{\frac{1}{2} HP_* (\mr{tr})} & 
HP_* (\mc O (\xnr (M) / W )) & \cong & H^{[*]} (\mh C^\times / (z \sim z^{-1})) \\
\downarrow & & \downarrow & & \downarrow \\
HP_* (\mc S_1^{\mf s}) & \longleftarrow & HP_* (C^\infty (\xunr (M)/W )) &
\cong & H^{[*]} (S^1 / (z \sim z^{-1} ))
\end{array}
\]
we see that the left vertical arrow is indeed an isomorphism.
From this we can derive a natural isomorphism
\[
HP_* (\mc H (G)^{\mf s}) \cong HP_* (\mc H (G,K)) \to 
HP_* (\mc S (G,K), \hot) \cong HP_* (\mc S (G)^{\mf s} ,\hot) \,.
\]
We remark that the Borel component is the most complicated
Bernstein component of $SL_2 (\mh Q_p )$. Indeed with one exception 
every other $\mf s \in \Omega (SL_2 (\mh Q_p ))$ has a trivial Weyl group, 
and therefore Prim$ \big( \mc H (SL_2 (\mh Q_p ))^{\mf s} \big)$ is 
homeomorphic to either $\mh C^\times$ or a point. Moreover both 
$\mc H (SL_2 (\mh Q_p ))^{\mf s}$ and $\mc S (SL_2 (\mh Q_p ))^{\mf s}$
are Morita equivalent to commutative algebras for such $\mf s$.
\vspace{4mm}

\section{Equivariant cosheaf homology}
\label{sec:3.3}

The most interesting applications of the results of Section \ref{sec:3.1} 
lie in their connection with the Baum--Connes conjecture. 
To make this relation precise we need several additional homology theories.
In particular, our forthcoming discussing will require some detailed knowledge
of equivariant cosheaf homology. Therefore we first provide an overview of this 
theory, which is mostly taken from \cite{BCH,HiNi}.

Let $G$ be a totally disconnected group and $\Sigma$ a polysimplicial
complex. We assume that $\Sigma$ is equipped with a polysimplicial $G$-action,
which is proper in the sense that the isotropy group $G_\sigma$ of any 
polysimplex $\sigma$ is compact and open. Let $\Sigma^p$ denote the collection 
of $p$-dimensional polysimplices of $\Sigma$, endowed with the discrete topology. 
Define the vector space 
\[
C_p (G ; \Sigma ) := {\ts \bigoplus_{\sigma \in \Sigma^p} } \, C_c^\infty (G_\sigma )
\]
where $C^\infty (X)$ denotes the set of locally constant complex valued functions
on a totally disconnected space $X$. We write the elements of $C_p (G;\Sigma )$
as formal sums $\sum_\sigma f_\sigma [\sigma ]$. If $\tau$ is a face of $\sigma$
then $G_\tau \supset G_\sigma$, so we may consider $f_\sigma$ as a locally constant
function on $G_\tau$. On every polysimplex we fix an orientation, and we identify
$[\bar \sigma]$ with $-[\sigma]$, where $\bar \sigma$ means $\sigma$ with the
opposite orientation. We write the simplicial boundary operator as
\[
\delta \sigma = {\ts \sum_{\tau \in \Sigma^{p-1}} } [\sigma : \tau] \, \tau \qquad 
\text{with} \qquad [\sigma : \tau ] \in \{ -1,0,1 \} \,.
\]
This gives a differential
\begin{align*}
& \delta_p : C_p (G ; \Sigma) \to C_{p-1} (G;\Sigma ) \,, \\
& \delta_p \big( f_\sigma [\sigma] \big) = {\ts \sum_{\tau \in \Sigma^{p-1}} }
[\sigma : \tau] f_\sigma [\tau] \,.
\end{align*}
We endow the differential complex $( C_* (G;\Sigma ) ,\delta_* )$ with
the $G$-action 
\[
g \cdot f_\sigma [\sigma ] = f_\sigma^g [g \sigma] \,,
\]
where $f_\sigma^g \in C_c^\infty (G_{g \sigma})$ is defined by
$f_\sigma^g (h) = f_\sigma (g^{-1} h g)$ and $g \sigma$ is endowed with the
orientation coming from our chosen orientation on $\sigma$. Clearly $\delta_*$ is 
$G$-equivariant, so it is well-defined on the space $C_* (G;\Sigma )_G$ of 
$G$-coinvariants. The equivariant cosheaf homology of $\Sigma$ is 
\begin{equation}
CH_n^G (\Sigma ) := H_n ( C_* (G;\Sigma )_G ,\delta_* ) \,.
\end{equation}
There is also a relative version of this theory. Let $\Sigma'$ be a $G$-stable
subcomplex of $\Sigma$. We define the relative equivariant cosheaf homology of
$(\Sigma ,\Sigma')$ as
\begin{equation}
CH_n^G (\Sigma ,\Sigma' ) := 
H_n ( C_* (G;\Sigma )_G / C_* (G;\Sigma')_G ,\delta_* ) \,.
\end{equation}
As usual there is a long exact sequence in homology:
\begin{equation}\label{eq:3.22}
\cdots \to CH_{n+1}^G (\Sigma ,\Sigma') \to CH_n^G (\Sigma' ) 
\to CH^G_n (\Sigma ) \to CH_n^G (\Sigma ,\Sigma') \to \cdots
\end{equation}
If $G$ acts freely on $\Sigma$ then $CH^G_n (\Sigma ,\Sigma')$ reduces
to the usual simplicial homology $H_n (\Sigma /G ,\Sigma' /G) $
with complex coefficients.

Higson and Nistor \cite{HiNi} introduced a natural map
\begin{equation}\label{eq:3.37}
CH_n (\Sigma ) \to HH_n (C_c^\infty (G)) \,,
\end{equation}
whose construction we recall in as much detail as we need. Let 
\[
\prefix{}{_n}{\hat \Sigma^p} := \{ (g_0 ,g_1 ,\ldots, g_n ,\sigma ) \in
G^{n+1} \times \Sigma^p : g_0 g_1 \cdots g_n \sigma = \sigma \}
\]
be the $n$th Brylinski space of $\Sigma^p$. By definition $\Sigma^p$ is
discrete, so $\prefix{}{_n}{\hat \Sigma^p}$ is a totally disconnected
space and $C_c^\infty \big( \prefix{}{_n}{\hat \Sigma^p} \big)$ is
defined. According to \cite[Section 4]{HiNi} there is an exact sequence
\begin{equation}\label{eq:3.15}
0 \leftarrow C_p (G;\Sigma )_G \leftarrow C_c^\infty \big( 
\prefix{}{_0}{\hat \Sigma^p} \big) \leftarrow C_c^\infty \big( 
\prefix{}{_1}{\hat \Sigma^p} \big) \leftarrow C_c^\infty \big(
\prefix{}{_2}{\hat \Sigma^p} \big) \leftarrow \cdots
\end{equation}
Consequently $CH_*^G (\Sigma )$ can be computed as the homology of a
double complex
\begin{equation}\label{eq:3.8}
\begin{array}{cccccc}
\downarrow & & \downarrow & &\downarrow \\
C_c^\infty \big( \prefix{}{_2}{\hat \Sigma^0} \big) & \leftarrow & C_c^\infty \big(
\prefix{}{_2}{\hat \Sigma^1} \big) & \leftarrow & C_c^\infty \big(
\prefix{}{_2}{\hat \Sigma^2} \big) &  \leftarrow \\
\downarrow & & \downarrow & &\downarrow \\
C_c^\infty \big( \prefix{}{_1}{\hat \Sigma^0} \big) & \leftarrow & C_c^\infty \big(
\prefix{}{_1}{\hat \Sigma^1} \big) & \leftarrow & C_c^\infty \big(
\prefix{}{_1}{\hat \Sigma^2} \big) &  \leftarrow \\
\downarrow & & \downarrow & &\downarrow \\
C_c^\infty \big( \prefix{}{_0}{\hat \Sigma^0} \big) & \leftarrow & C_c^\infty \big(
\prefix{}{_0}{\hat \Sigma^1} \big) & \leftarrow & C_c^\infty \big(
\prefix{}{_0}{\hat \Sigma^2} \big) &  \leftarrow 
\end{array}
\end{equation}
In this diagram the horizontal maps come from the boundary map $\partial$ on
$\Sigma$, while the vertical maps are essentially the differentials for a 
Hochschild complex. There are natural maps 
\begin{equation}\label{eq:3.16}
\begin{array}{lll}
C_c^\infty \big( \prefix{}{_n}{\hat \Sigma^0} \big) & \to & C_c^\infty (G^{n+1}) \\
f & \to & \big( (g_0 ,\ldots ,g_n) \mapsto \sum_{x \in \Sigma^0} f (g_0 ,\ldots ,g_n, x) \big)
\end{array}
\end{equation}
from the double complex \eqref{eq:3.8} to the standard Hochschild complex for
$C_c^\infty (G)$. Together \eqref{eq:3.15}, \eqref{eq:3.8} and \eqref{eq:3.16} 
yield the map \eqref{eq:3.37}.

Suppose now that $G$ acts properly on some affine building. The Hochschild homology 
of $C_c^\infty (G)$ admits a decomposition 
\[
HH_n (C_c^\infty (G)) = HH_n (C_c^\infty (G))_{\mr{ell}} \oplus HH_n (C_c^\infty (G))_{\mr{hyp}}
\]
into an elliptic and an hyperbolic part. Upon periodization the hyperbolic part
disappears and one finds \cite[Section 7]{HiNi}
\begin{equation}\label{eq:3.17}
HP_n (C_c^\infty (G)) = {\ts \bigoplus_{m \in \mh Z} } HH_{n + 2m} (C_c^\infty (G))_{\mr{ell}} \,.
\end{equation}
This and \eqref{eq:3.37} yield a map 
\begin{equation}\label{eq:3.12}
\mu_{HN} : CH_{[*]}^G (\Sigma ) \to HP_* (C_c^\infty (G)) \,.
\end{equation}
Let $h : X \to \Sigma$ be any morphism of proper $G$-simplicial complexes. 
From the explicit formula \eqref{eq:3.16} we see that $\mu_{HN}$ for $X$ factors as
\begin{equation}\label{eq:3.38}
CH_n^G (X) \xrightarrow{CH_n^G (h)} CH_n^G (\Sigma ) \to HH_n (C_c^\infty (G)) 
\to HP_n (C_c^\infty (G))\,.
\end{equation}
Higson and Nistor \cite{HiNi} showed that \eqref{eq:3.12} is an isomorphism if 
$\Sigma$ is an affine building. The special case where $G$ is a simple $p$-adic 
group was also proved by Schneider \cite{Schn}. 
\vspace{4mm}

\section{The Baum--Connes conjecture}
\label{sec:3.4}

Let $G$ be any locally compact group acting properly on a Hausdorff space $\Sigma$. 
A subspace $X \subset \Sigma$ is called $G$-compact if $X / G$ is compact. 
The equivariant $K$-homology of $\Sigma$ is defined as
\begin{equation}\label{eq:3.18}
K_*^G (\Sigma ) = \varinjlim KK_*^G (C_0 (X), \mh C )
\end{equation}
where $KK_*^G$ is Kasparov's equivariant $KK$-theory \cite{Kas}
and the limit runs over all $G$-compact subspaces $X$ of $\Sigma$.
The Baum--Connes conjecture asserts that the assembly map
\begin{equation}\label{eq:3.7}
\mu : K_*^G (\Sigma) \to K_* (C_r^* (G))
\end{equation}
is an isomorphism if $\Sigma$ is a classifying space for proper $G$-actions.
Building upon the work of Kasparov, Vincent Lafforgue proved this conjecture 
for many groups, including all locally compact groups 
that act properly isometrically on an affine building \cite{Laf}.

Now we specialize to a reductive $p$-adic group $G$. In this case the affine 
Bruhat--Tits building $\beta G$ is a classifying space for proper $G$-actions. 
We recall that $\beta G$ is a finite dimensional locally finite polysimplicial 
complex endowed with an isometric $G$-action such that $\beta G / G$ 
is compact and contractible.
For any $X$ as above there exists a continuous $G$-map $h : X \to \beta G$,
and it is unique up to homotopy. The assembly map 
$\mu : K_*^G (X) \to K_* (C_r^* (G)$ factors as
\begin{equation}\label{eq:3.36}
K_*^G (X) \xrightarrow{K_*^G (h)} K_*^G (\beta G) \xrightarrow{\mu} K_* (C_r^* (G)) \,.
\end{equation}
A natural receptacle for a Chern character from $K_*^G (X)$ is formed by 
\[
HL_*^G (X) := HL_*^G (C_0 (X), \mh C)
\]
where $HL_*^G$ denotes equivariant local cyclic homology, as defined
and studied by Voigt \cite{Voi3}.
With these notions we can state and prove a more precise
version of \cite[Proposition 9.4]{BHP2}. We note that a similar idea
was already used in \cite{BHP1} to prove the Baum--Connes 
conjecture for $G = GL_n (\mh F )$.

\begin{thm}\label{thm:3.4}
There exists a commutative diagram 
\[
\begin{array}{ccc}
K_*^G (\beta G) & \xrightarrow{\quad \mu \quad} & K_* (C_r^* (G)) \\
\downarrow ch & & \downarrow ch \\
HL_*^G (\beta G) & & HP_* (\mc S (G),\hot_{\mh C} ) \\
\downarrow & & \uparrow \\
CH_{[*]}^G (\beta G) & \xrightarrow{\quad \mu' \quad} & HP_* (\mc H (G)) 
\end{array}
\]
with the properties:
\begin{description}
\item[a)] Both Chern characters become isomorphisms after 
applying $\otimes_{\mh Z} \mh C$ to their domain.
\item[b)] The other maps are natural isomorphisms.
\end{description}
\end{thm}
\emph{Proof.} As mentioned before, Lafforgue \cite{Laf} showed that the 
assembly map $\mu$ is an isomorphism. The right column is taken 
care of by Theorems \ref{thm:3.2} and \ref{thm:3.3}. 

Let $\Sigma$ be a finite dimensional, locally finite $G$-compact 
$G$-simplicial complex.
It was proved in \cite[Proposition 10.4]{Voi3} that the inclusion map
$C_c^\infty (\Sigma) \to C_0 (\Sigma)$ induces an isomorphism
\[
HL^G_* (\Sigma) =  HL_*^G (C_0 (\Sigma),\mh C ) \to
HL_*^G (C_c^\infty (\Sigma),\mh C ) \,.
\]
Let $HP_*^G$ denote Voigt's equivariant periodic cyclic homology \cite[Section 3]{Voi2}. 
According to \cite[Section 6]{Voi4} there are natural isomorphisms
\[
\begin{array}{rcl}
HL_*^G (C_c^\infty (\Sigma),\mh C ) & \cong & HP_*^G (C_c^\infty (\Sigma),\mh C ) \,, \\
ch \,:\, KK_*^G (C_0 (\Sigma) ,\mh C ) \otimes_{\mh Z} \mh C & \to & 
 HL_*^G (C_0 (\Sigma) ,\mh C ) \,.
\end{array}
\]
We have to check that
\begin{equation}\label{eq:3.9}
HP_*^G (C_c^\infty (\Sigma) ,\mh C ) \cong CH_{[*]}^G (\Sigma) \,.
\end{equation}
Baum and Schneider \cite[Section 1.B]{BaSc} showed that cosheaf homology 
can be regarded as a special case of a bivariant (co)homology theory:
\begin{equation}\label{eq:3.10}
CH_n^G (\Sigma) \cong H^n_G (\Sigma, \mr{point}) \,.
\end{equation}
According to \cite{Voi2} the right hand side of \eqref{eq:3.10} is naturally
isomorphic to the left hand side of \eqref{eq:3.9}, so we get natural isomorphisms
\begin{equation}\label{eq:3.31}
K_*^G (\Sigma ) \otimes_{\mh Z} \mh C \to HL_*^G (\Sigma ) \to CH_{[*]}^G (\Sigma ) \,.
\end{equation}
The case $\Sigma = \beta G$ gives us the left column of the theorem. 
To complete the proof we define 
\begin{equation}\label{eq:3.13}
\mu' : CH^G_{[*]} (\beta G) \to HP_* (\mc H (G))
\end{equation}
as the unique map so that the diagram commutes. $\qquad \Box$
\\[2mm]

It is not immediately clear that \eqref{eq:3.13} and \eqref{eq:3.12}, for
$C_c^\infty (G) = \mc H (G)$ and $\Sigma = \beta G$, are the same map. 
We will prove this by reduction to the following simpler case.
Let $U \in$ CO$ (G)$ and consider the discrete proper homogeneous
$G$-space  $G/U$. By the universal property of $\beta G$ there exists
a continuous $G$-equivariant map $G/U \to \beta G$, and it is unique
up to homotopy. With a suitable simplicial subdivision of $\beta G$ we can
achieve that this is in fact a simplicial $G$-map.

\begin{lem}\label{lem:3.5} 
The following diagram commutes for elements in the upper left corner.
\[
\begin{array}{ccccc}
K_*^G (G/U) & \to & K_*^G (\beta G) & \xrightarrow{\mu} & K_* (C_r^* (G)) \\
\downarrow & &\downarrow & &\downarrow \\
CH^G_{[*]} (G/U) & \to & CH^G_{[*]} (\beta G) & 
\xrightarrow{\mu_{HN}} & HP_* (\mc H (G))
\end{array}
\]
\end{lem}
\emph{Proof.} 
The left hand square commutes by functoriality. It follows readily from the 
definitions that $K_*^G (G/U) \cong K_*^U (\mr{point})$. By the 
functoriality of the Baum-Connes assembly map there is a commutative diagram
\begin{equation}\label{eq:3.20}
\begin{array}{ccc}
K_*^U (\mr{point}) &\xrightarrow{\mu_U} & K_* (C_r^* (U)) \\
\downarrow & & \downarrow \\
K_*^G (\beta G) & \xrightarrow{\mu} & K_* (C_r^* (G)) \,.
\end{array}
\end{equation}
Since $U$ is compact and totally disconnected, both $K_*^U (\mr{point})$ 
and $K_* (C_r^* (U))$ are naturally isomorphic to the ring of smooth (virtual) 
representations $R (U)$, and $\mu_U$ corresponds to the composition 
of these isomorphisms. The right vertical map comes from the inclusion 
$C_r^* (U) \to C_r^* (G)$, so it sends a $U$-module $V$ to $\mr{Ind}_U^G (V)$.\\ 
Similarly there are a canonical isomorphism 
\[
CH_*^G (G/U) \cong CH_*^U (\mr{point})
\]
and a commutative diagram
\begin{equation}\label{eq:3.21}
\begin{array}{ccc}
CH_{[*]}^U (\mr{point}) &\xrightarrow{\mu_{HN}} & HP_* (C_c^\infty (U)) \\
\downarrow & & \downarrow \\
CH_{[*]}^G (\beta G) & \xrightarrow{\mu_{HN}} & HP_* (C_c^\infty (G)) \,.
\end{array}
\end{equation}
According to \cite[Section 4]{HiNi} we have $HP_1 (C_c^\infty (U)) = 0$ and
\[
HP_0 (C_c^\infty (U)) = HH_0 (C_c^\infty (U)) = C_c^\infty (U)_U \,.
\]
By definition, also
\[
CH_n^U (\mr{point}) = \left\{ \begin{array}{lll}
C_c^\infty (U)_U & \mr{if} & n = 0 \\
0 & \mr{if} & n > 0 \,.
\end{array} \right.
\]
A glance at the double complex \eqref{eq:3.8} shows that 
\[
\mu_{HN} : CH_0^U (\mr{point}) \to HP_0 (C_c^\infty (U))
\]
corresponds to the identity map under these identifications.
Furthermore $U$ is profinite, so
\[
C_c^\infty (U) = \varinjlim \mh C [F] ,
\]
where the limit runs over all finite quotient groups $F$ of $U$. Similarly we can write 
$C_r^* (U)$ as an inductive limit in the category of $C^*$-algebras. In this situation
both $K_*$ and $HP_*$ commute with $\varinjlim$, so we get a Chern character
\begin{equation}\label{eq:3.33}
K_* (C_r^* (U)) \cong \varinjlim K_* (\mh C [F]) \to 
 \varinjlim HP_* (\mh C [F]) \cong HP_* (C_c^\infty (U)) \,.
\end{equation}
Now the right hand square of the diagram
\begin{equation}\label{eq:3.19}
\begin{array}{ccccc}
K_*^U (\mr{point}) & \xrightarrow{\mu} & K_* (C_r^* (U)) &
\to & K_* (C_r^* (G)) \\
\downarrow & & \downarrow & & \downarrow \\
CH_{[*]}^U (\mr{point}) & \xrightarrow{\mu_{HN}} & HP_* (C_c^\infty (U)) &
\to & HP_* (\mc H (G))
\end{array}
\end{equation}
commutes by functoriality.
According to Voigt \cite[Proposition 13.5]{Voi3} the Chern character
\[
K_*^U (\mr{point}) \to HL_*^U (\mr{point})
\]
can be identified with the character map $R (U) \to C_c^\infty (U)^U$.
The isomorphism between $HL_*^U$(point) and $CH_*^U$(point)
then becomes the canonical map
\[
C_c^\infty (U)^U \to C_c^\infty (U)_U ,
\]
which is bijective because $U$ is compact. Since the Chern character for $\mh C [F]$ 
in \eqref{eq:3.33} may also be identified with the character map, we find that the 
left hand square of \eqref{eq:3.19} commutes. Together the commutative diagrams 
\eqref{eq:3.20}, \eqref{eq:3.21} and \eqref{eq:3.19} complete the proof. $\qquad \Box$
\\[2mm]
\begin{lem}\label{lem:3.6}
The maps $\mu_{HN}$ from \eqref{eq:3.12} and $\mu'$ from \eqref{eq:3.13} are the same.
\end{lem}
\emph{Proof.}
We have to show that the diagram
\begin{equation}\label{eq:3.14}
\begin{array}{ccc}
K_*^G (\Sigma) \otimes_{\mh Z} \mh C & \xrightarrow{\mu} & 
K_* (C_r^* (G)) \otimes_{\mh Z} \mh C \\
\downarrow & &\downarrow \\
CH^G_{[*]} (\Sigma) & \xrightarrow{\mu_{HN}} & HP_* (\mc H (G))
\end{array}
\end{equation}
commutes for $\Sigma = \beta G$.
By subdividing all polysimplices we may assume that $\beta G$ is a 
simplicial complex, and that $G$ preserves this structure. Let $\beta^{(n)} G$ 
denote the $n$-skeleton of $\beta G$, and $(\beta G)^n$ the collection of $n$-simplices. 
Both inherit a $G$-action from $\beta G$.
We will prove the commutativity of \eqref{eq:3.14} for $\Sigma = \beta^{(n)}$, with 
induction to $n$. 

The set $\beta^{(0)} G$ is a finite union of $G$-spaces of the form $G/U$ with 
$U \in$ CO$ (G)$, so the case $n = 0$ follows from Lemma \ref{lem:3.5}.
Similarly \eqref{eq:3.14} commutes for $\Sigma = (\beta G)^n$.

We consider $(\beta G)^n \times S^n$ as a $G$-space with a trivial action on $S^n$. The
long exact sequence \eqref{eq:3.22} for the pair 
$\big( (\beta G)^n \times S^n , (\beta G)^n \times \mr{point} \big)$ reads
\[
\begin{split}
\cdots & \to CH_p^G ((\beta G)^n ) \to CH_p^G ((\beta G)^n \times S^n ) \xrightarrow{\phi}
CH_{p-n}^G ((\beta G)^n ) \otimes_{\mh C} H_n (S^n ,\mr{point}) \\
& \to CH_{p-1}^G ((\beta G)^n ) \to CH_{p-1}^G ((\beta G)^n \times S^n ) \to \cdots
\end{split}
\]
We note that \pagebreak[3]
\begin{equation*}
CH_p^G ((\beta G)^n \times S^n ) \cong CH_p^G ((\beta G)^n ) \otimes_{\mh C} H_0 (S^n ) 
\oplus CH_{p-n}^G ((\beta G)^n ) \otimes_{\mh C} H_n (S^n ,\mr{point}) ,
\end{equation*}
so $\phi$ is surjective and the sequence splits. Since $\mu_{HN} - \mu' = 0$ on 
$CH_*^G ((\beta G)^n )$, there exists a unique map $f$ making the following diagram commutative:
\begin{equation}\label{eq:3.34}
\begin{array}{ccccc}
\!\! CH_p^G ((\beta G)^n ) \!\! & \to & \!\! CH_p^G ((\beta G)^n \times S^n ) \!\! & \to &
\!\! CH_{p-n}^G ((\beta G)^n ) \! \otimes_{\mh C} \! H_n (S^n ,\mr{point}) \!\! \\
\downarrow 0 & & \downarrow \mu_{HN} - \mu' & & \downarrow f \\
HP_* (\mc H (G)) & = & HP_* (\mc H (G)) & = & HP_* (\mc H (G)) \,. \\
\end{array}
\end{equation}
By the universal property of $\beta G$, there exists a $G$-map 
$h : (\beta G )^n \times S^n \to \beta G$, and it is unique up to homotopy. Hence we may 
assume that $h$ maps $\{ \sigma \} \times S^n \subset (\beta G )^n \times S^n$ to the 
barycenter of $\sigma$ in $\beta G$. By \eqref{eq:3.38} and \eqref{eq:3.36} the middle map 
in \eqref{eq:3.34} factors as
\begin{equation}\label{eq:3.35} 
CH_p^G ((\beta G)^n \times S^n ) \xrightarrow{CH_p^G (h)} CH_p^G (\beta G )  
\xrightarrow{\mu_{HN} - \mu'} HP_* (\mc H (G)) \,.
\end{equation}
But $CH_p^G (h)$ kills the $n$th homology of $S^n$, which in combination with \eqref{eq:3.34}
shows that \eqref{eq:3.35} is zero. Therefore the above map $f$ must also be zero.
As $G$-spaces we have 
\[
\beta^{(n)} G \setminus \beta^{(n-1)} G \cong 
(\beta G)^n \times S^n \setminus \mr{point} \,,
\]
with $G$ acting trivially on the last factor. The corresponding long exact sequence in
equivariant cosheaf homology is
\[
\begin{split}
\cdots & \to CH_p^G ( \beta^{(n-1)} G ) \to CH_p^G ( \beta^{(n)} G ) \to
CH_{p-n}^G ((\beta G)^n ) \otimes_{\mh C} H_n (S^n ,\mr{point}) \\
& \to CH_{p-1}^G (\beta^{(n-1)} G ) \to CH_{p-1}^G ( \beta^{(n)} G ) \to \cdots
\end{split}
\]
By the induction hypothesis $\mu_{HN} - \mu' = 0$ on $CH_*^G ( \beta^{(n-1)} G )$, so
we can write down a commutative diagram
\[
\begin{array}{ccccc}
CH_p^G ( \beta^{(n-1)} G) & \to & CH_p^G (\beta^{(n)} G) & \to &
CH_{p-n}^G ((\beta G)^n ) \otimes_{\mh C} H_n (S^n ,\mr{point}) \\
\downarrow 0 & & \downarrow \mu_{HN} - \mu' & & \downarrow f \\
HP_* (\mc H (G)) & = & HP_* (\mc H (G)) & = & HP_* (\mc H (G)) \,. \\
\end{array}
\]
We already showed that $f = 0$, so $\mu_{HN} - \mu' = 0$ on $CH_p^G (\beta^{(n)} G)$.
Thus \eqref{eq:3.14} commutes for $\Sigma = \beta^{(n)} G$, which completes our induction 
step. $\qquad \Box$
\\[2mm]

The above proof can be compared with \cite[Section 5.1]{Mey1}.

\begin{cor}\label{cor:3.7}
It can be proved with periodic cyclic homology that the Baum--Connes assembly map
\[
\mu \otimes \mr{id} : K_*^G (\beta G) \otimes_{\mh Z} \mh Q \to 
K_* (C_r^* (G)) \otimes_{\mh Z} \mh Q
\]
is an isomorphism for every reductive $p$-adic group $G$.
\end{cor}
\emph{Proof.}
Lemma \ref{lem:3.6} and \eqref{eq:3.12} show that $\mu' = \mu_{HN}$ is an 
isomorphism. Hence all the isomorphisms in the diagram of Theorem 
\ref{thm:3.4} admit mutually independent proofs. With the commutativity of
the diagram we can use any five of them to prove the sixth. 
In particular we can show without using Lafforgue's work that 
\[
K^G_* (\beta G) \otimes_{\mh Z} \mh C \cong 
K_* (C_r^* (G)) \otimes_{\mh Z} \mh C \,,
\]
which is equivalent to $\mu$ being a rational isomorphism. $\qquad \Box$

\end{document}